\documentclass[10pt,reqno]{amsart} 
\usepackage[T1]{fontenc}
\usepackage{amssymb,alltt}
\usepackage[all]{xy}            
\usepackage{tikz}
    \usetikzlibrary{shadows,matrix,arrows}
    \newcommand*\ccircled[1]{\tikz[baseline=(C.base)]\node[draw,rectangle,rounded corners,inner sep=1.9pt,line width=0.15mm](C){\normalsize#1};\!} 
\usepackage{hyperref} 
	\hypersetup{pdflinkmargin=0pt}
\usepackage{enumitem}
\usepackage{multirow} 
\usepackage{verbatim} 
\usepackage{stackrel}
\newtheorem{Thm}[subsection]{Theorem}

\newtheorem{Lmm}[subsection]{Lemma}
\numberwithin{equation}{section}

\def\id{\mathrm{id}}        
\def\Ker{\mathrm{Ker}}      
\def\Im{\mathrm{Im}}      	
\def\Coker{\mathrm{Coker}}  
\def\N{\mathbb{N}}
\def\Z{\mathbb{Z}}
\def\Q{\mathbb{Q}}
\def\R{\mathbb{R}}
\def\C{\mathbb{C}}
\newcommand{\rd}[1]{\textcolor[rgb]{1.00,0.00,0.00}{#1}}
\newcommand{\gr}[1]{\textcolor[rgb]{0.6,0.6,0.6}{#1}}
\newcommand{\mns}{\!\raisebox{-0.5pt}{\scalebox{0.5}[1.0]{\( - \)}}\!}

\newcommand{\leftrarrows}{\mathrel{\raise.75ex\hbox{\oalign{ $\scriptstyle\leftarrow$\cr	\vrule width0pt height.5ex$\hfil\scriptstyle\relbar$\cr}}}}
\newcommand{\lrightarrows}{\mathrel{\raise.75ex\hbox{\oalign{ $\scriptstyle\relbar$\hfil\cr	$\scriptstyle\vrule width0pt height.5ex\smash\rightarrow$\cr}}}}
\newcommand{\Rrelbar}{\mathrel{\raise.75ex\hbox{\oalign{ $\scriptstyle\relbar$\cr \vrule width0pt height.5ex$\scriptstyle\relbar$}}}}
\newcommand{\longleftrightarrows}{\leftrarrows\joinrel\Rrelbar\joinrel\lrightarrows}
\makeatletter\def\leftrightarrowsfill@{\arrowfill@\leftrarrows\Rrelbar\lrightarrows}
\newcommand{\xleftrightarrows}[2][]{\ext@arrow 3399\leftrightarrowsfill@{#1}{#2}}\makeatother

\begin{document}
\title[Chain complex reduction via fast digraph traversal]{Chain complex reduction\\ via fast digraph traversal}
\author{Leon Lampret}
\address{Faculty of Mathematics and Physics, Department of Mathematics, University of Ljubljana, Jadranska ulica 21, 1000 Ljubljana, Slovenia}
\email{lampretl@gmail.com, leon.lampret@fmf.uni-lj.si, +386 1 47 66 624}
\date{October 1, 2020}

\begin{abstract}
Reducing a chain complex (whilst preserving its homotopy-type) using algebraic Morse theory gives the same end-result as Gaussian elimination, but AMT does it only \emph{on certain rows/columns} and with \emph{several pivots} (in all matrices simultaneously). Crucially, instead of doing costly row/column operations on a sparse matrix, it computes traversals of a bipartite digraph. This significantly reduces the running time and memory load (smaller fill-in and coefficient growth of the matrices). However, computing with AMT requires the construction of a valid set of pivots (called a Morse matching).
\par We discover a family of Morse matchings on \emph{any} chain complex of free modules of finite rank.
We show that \emph{every} acyclic matching is a \emph{subset} of some member of our family, so all maximal Morse matchings are of this type. 
\par Both the input and output of AMT are chain complexes, so the procedure can be used \emph{iteratively}. When working over a field or a local PID, this process ends in a chain complex with zero matrices, which produces homology. However, even over more general rings, the process often reveals homology, or at least reduces the complex so much that other algorithms can finish the job. Moreover, it also returns homotopy equivalences to the reduced complexes, which reveal the \emph{generators} of the homology and the \emph{induced maps} $H_\ast(\varphi)$.
\par We design a new algorithm for reducing a chain complex and implement it in \textsc{Mathematica}. We test that it outperforms other CASs. As a special case, given a sparse matrix over any field, the algorithm offers a new way of computing the rank and a sparse basis of the kernel (or null space), cokernel (or quotient space, or complementary subspace), image, preimage, sum and intersection subspace. It outperforms built-in algorithms in other CASs.
\end{abstract}

\maketitle

\subsection*{Motivation} A \textit{chain complex} is a sequence $C_\ast\!:C_0\!\overset{\partial_1}{\leftarrow}\! C_1\!\overset{\partial_2}{\leftarrow}\! C_2\!\overset{\partial_3}{\leftarrow}\!\ldots$ of $R$-module morphisms with the property that $\partial_k\!\circ\!\partial_{k+1}=0$ for all $k$; it is the central object of study in homological algebra. This branch of mathematics is developing ways to associate to a mathematical object $X$ a chain complex $C_\ast(X)$ (in a functorial way), and then compute its \textit{homology} modules $H_kC_\ast= \frac{\Ker\,\partial_k}{\Im\,\partial_{k+1}}$. These modules give a lot of valuable information about $X$.
\par (Co)homology theories were developed for objects $X$ coming from many diverse areas (e.g. a CW complex, topological space, smooth manifold, simplicial filtration, group, associative algebra, Lie algebra, knot or link, sheaf, algebraic set, and many more). Hence, there is a big demand for efficient methods (time and memory wise) to calculate homology. Unfortunately, this continues to be a hard problem, as it is still infeasible to compute with sparse matrices $\partial_k$ larger than $10^6\!\times\!10^6$.
\par Let us see that in general, \textit{the bigger the matrices in $C_\ast$ get, the emptier they become}. Take for instance the Poincar\'{e} \cite[pp.4]{bookWeibelHA} / Eilenberg-MacLane \cite[pp.177]{bookWeibelHA} / Hochschild \cite[pp.300]{bookWeibelHA} / Chevalley \cite[pp.238]{bookWeibelHA} chain complex with trivial coefficients. The density $\rho=\frac{\text{number of nonzero entries}}{\text{number of all entries}}$ of the $k$-th matrix is:
\begin{itemize} \setlength{\itemindent}{-8mm}
	\item \!\!$\rho(\partial_k)\!=\!\frac{(k+1)\cdot f_k}{f_{k\!-\!1}\cdot f_k}\xrightarrow{f_{k\!-\!1}\to\infty}0$, where $f_k$ = number of $k$-faces of a simplicial complex;
	\item $\rho(\partial_k)\!=\!\frac{(k+1)\cdot n^k}{n^{k\!-\!1}\cdot n^k}\xrightarrow{n\to\infty}0$, where $n$ = cardinality of a group;
	\item $\rho(\partial_k)\!<\!\frac{(k+1)n\cdot n^{k+1}}{n^k\cdot n^{k+1}}\xrightarrow{n\to\infty}0$, where $n$ = dimension of an associative algebra;
	\item $\rho(\partial_k)\!<\!\frac{\binom{k}{2}n\cdot \binom{n}{k}}{\binom{n}{k\!-\!1}\cdot \binom{n}{k}}\xrightarrow{n\to\infty}0$, where $n$ = dimension of a Lie algebra.
\end{itemize}
Lesson learned: it makes sense to focus on methods for \emph{sparse} chain complexes.

\subsection*{Prior work} In the last 40 years, a lot of attention was devoted to this area, in theoretical as well as applied circles. We list chronologically a few popular methods.
\par If $R$ is a field, then the homology is computed from the ranks of sparse matrices: $\dim H_kC_\ast=\dim C_k\!-\!\mathrm{rank}\,\partial_k\!-\!\mathrm{rank}\,\partial_{k+1}$. Efficient computation of rank is explored in \cite{articleWiedemannSSLEFF, articleMulmuleyFPACRMAF, articleDumasVillardCRLSMFF, articleDuranSaundersWanHARSM, articleDumasVincentGiorgiUrbanskaPCRLSMAK, articleMaySaundersWanEMRCASSRG, articleGotsmanToledoCNSSRM, articleSaundersYouseLMSR, articleCheungKwokLauFMRAA, articleFosterDavisA933, articleBouillaguetDelaplaceSGEMpU, articleSaadUbaruFMENRLM, articleBouillaguetDelaplaceVogePSpluqFMp, articleGondzioSchorkRRGEMVC}.
If $R$ is a subfield of $\C$, then homology equals the number of zero eigenvalues of the Laplacian operators: $\dim H_kC_\ast=\dim C_k-\mathrm{rank}(\partial_k^h\!\circ\!\partial_k\!+\!\partial_{k+1}\!\circ\!\partial_{k+1}^h)$, where $^h$ is the conjugate transpose of a matrix. This method was developed in \cite{articleEckmannHFREK, articleKostantLACGBWT, articleHozoELHLACP, articleSiggLHFTNLA, articleFriedmanCBNCL, articleKookReinerStantonCLMC, articleDuvalReinerSSCLI, articleHanlonLM, articleMaleticRajkovicCLESCACN, articleHorakJostSCLOSC}.
\par If $R\!=\!K[t]$ where $K$ is a field and $C_k$ are graded modules (this is the case for persistent homology of filtered chain complexes), then there is a plethora of highly efficient but specialized methods of computing homology. The most important ones are described in \cite{articleCarlssonZomorodianCPH, articleSilvaMorozovJohanssonDPC, articleMischaikowNandaMTFECPH, articleAdamsTauszJohanssonJP, articleBinchiMerelliRuccoPetriVaccarinoJH, articleMariaBoissonnatGlisseYvinecGL, articleBauerKerberReininghausWagnerPHAT, articleOtterPorterTillmannGrindrodHarringtonRCPH, articleBauer, articleLutgehetmannGovcSmithLeviCPHDFC}.
\par If $R$ is a PID (previous paragraphs are special cases of this), then the homology is computed from the Smith normal forms: $H_kC_\ast\cong R^r\!\oplus\!\bigoplus_{t\in T}\!\frac{R}{tR}$, where $r$ is the number of zero diagonals in the SNF of $\partial_k$ minus the number of nonzero diagonals in $\partial_{k+1}$ and $T$ is the multiset of diagonal entries in the SNF of $\partial_{k+1}$ that are not $0$ or a unit in $R$. Efficient computation of SNFs was researched in \cite{articleBarnettPaceEALSC1, articleBachemKannanPACSHNFIM, articleSmithCSNFIM, articleKaltofenKrishnamoorthySaundersFPCHSFPM, articleKaltofenKrishnamoorthySaundersSV, articleIliopoulosWCBACCSFAGHSNFIM, articleKaltofenKrishnamoorthySaundersPAMNF, articleVillardCSNFPM, articleVillardFPCSNFPM, articleVillardGSCSNFPM, articleStorjohannNOACSNFIM, articleStorjohannLabahnFVACSNFPM, articleVillardFPAMRNF, articleDumasSaundersVillardISMV, articleEberlyGiesbrechtVillardCDSNFIM, articleDumasSaundersVillardESIMSNFC, articleGiesbrechtFCSFSIM, articleDumasHeckenbachSaundersCSHBESNFA, articleJagerPACSNFLM, articleJagerNACSNFIPM, articleVillardSRPELARQ, articleJagerWagnerEPHSNFA, articleWilkeningLCSNFMP, articleElsheikhGiesbrechtNovocinSaundersFCSNFSMLR}.
\par In more general situations, e.g. if $R\!=\!K[x_1,\ldots,x_n]/(f_1,\ldots,f_m)$ is a quotient polynomial ring and $C_k$ are finitely generated, homology can be computed using Gr\"{o}bner bases \cite[p.612, \texttt{homolog.lib}]{bookGreuelPfisterSICA}, but since there is currently no classification of modules over such rings, those results are less relevant for our goals.
\par Special cases of the content of this article and other related work has also been done in \cite{articleDlotkoWagnerCHPHIMD, articleKaczynskiMrozekSlusarekHCRCC, articleDiazRealCCOFSC, articlePilarczykRealCCHCCOCC, 
	articleBatkoMrozekCHA, articleDlotkoKaczynskiMrozekWannerCHARC, articleMrozekPilarczykZelaznaHAAS, 	 articleCurryGhristNandaDMTCCSC, articleHarkerMischaikowMrozekNandaDMTACHCM, articleBauerRathodHAMM}.

\subsection*{Results} All previous approaches to computing (co)homology use matrices. But row/column operations on a sparse matrix are costly, since they reconfigure the whole data structure every time. We dispense with that strategy and work using bipartite acyclic digraphs with weighted edges. All that is needed is to compute, for a small subset of start-point vertices, the end-points and weights of paths. The digraph data structure remains \emph{unaltered} during the computation. 
\par Using AMT, we design a new algorithm \ref{5.1.algorithm}, which does not deal with boundary matrices one at a time, but reduces the size of all $\partial_k$ \emph{simultaneously} (while still preserving the homotopy type and hence homology of the chain complex). In this process, $\partial_{k-1}$ and $\partial_{k+1}$ help in the reduction of $\partial_k$. For every invertible entry $\partial_{k,i,j}$ which has zeros left and below it (i.e. $\partial_{k,i,<j}\!=\!0\!=\!\partial_{k,>i,j}$), AMT deletes the $i$-row of $\partial_k$, $i$-column of $\partial_{k-1}$, $j$-column of $\partial_k$, $j$-row of $\partial_{k+1}$, and (usually) alters the remaining entries in $\partial_k$. Before the computation even begins, many rows and columns in the boundary matrices can have their entries set to zero, without affecting the outcome, thereby saving up time and memory. 
\par The capabilities of our implementation are documented in \ref{5.7.performance}. It requires no computation of ranks or SNFs of matrices (though they can be reconstructed from our results). Several examples from simplicial complexes as well as associative and Lie algebras show that we can efficiently calculate (on a laptop computer) the homology of complexes, which contain matrices larger than $10^7\!\times\!10^7$. One of the advantages of AMT is that the input and output of the algorithm are both chain complexes, hence our procedure can be applied \emph{iteratively} (unlike smooth / discrete Morse theory, which transforms a smooth manifold / simplicial complex into a CW complex with semi-known gluing maps). 
\par All these calculations are not possible without Morse matchings. Before, these were constructed either algorithmically (e.g. \cite[p.887]{articleCurryGhristNandaDMTCCSC}, \cite[p.165]{articleHarkerMischaikowMrozekNandaDMTACHCM}, \cite{articleBauerRathodHAMM}, \cite[p.13]{articleDlotkoWagnerCHPHIMD}) or explicitly but case by case (e.g. \cite[p.123]{articleSkoldbergMTFAV}, \cite[p.44,53,68]{articleJollenbeckADMTACA} \cite[p.265]{articleLampretVavpeticCLAAMT}, \cite[p.8,10,15]{articleLampretVavpeticCPLA}, \cite[p.2]{articleLampretVavpeticHCEAAMT}). Now, in \ref{2.1.formulation}, we discover an explicit universal way of constructing Morse matchings. We prove that AMT is a \emph{functor} and show how it computes (co)homology with generators and induced maps, the torsion in homology, and (co)homology operations. The generality of our construction offers theoretical advances for appropriate homology theories, which will be shown in future papers.

\subsection*{Conventions} Throughout this article, $R$ will be a commutative unital ring and\vspace{-1mm} \[C_\ast\!: C_0\!\overset{\partial_1}{\leftarrow}\!C_1\!\leftarrow\!\ldots\! \leftarrow\!C_{N\!-\!1} \!\overset{\partial_N}{\leftarrow} \!C_N\vspace{-1mm}\] a chain complex of free $R$-modules of finite rank. Also, $R^\times$ is the group of units (=invertible elements) of $R$, e.g. $\Z^\times\!=\!\{1,-1\}$ and $K^\times\!=\!K\!\setminus\!\{0\}$ for any field $K$. 
\vspace{4mm}

\section{Formulation of AMT}
\noindent Our algorithm uses algebraic Morse theory (which is a generalization of Forman's discrete Morse theory \cite{articleFormanMTCC}), so we include a concise review of it. The exposition is rather technical, but the invested diligence will pay off.
\subsection{Formulation}\label{formulation} Pick a basis $I_k$ for each $C_k$. Let $\Gamma_{\!C_\ast}$ be a graph, with vertex set the disjoint union $\bigsqcup_kI_k$, and for every nonzero entry $\partial_{k,u,v}\!=\!w \in R\!\setminus\!\{0\}$ in the matrix $\partial_k$ a directed weighted edge $u\!\overset{w}{\leftarrow}\!v$.
\par A \emph{Morse matching} is any collection $\mathcal{M}$ of edges from $\Gamma_{C_\ast}$, such that:
\begin{enumerate}
	\item[$(1)$] $\mathcal{M}$ is a \emph{matching}, i.e. edges in $\mathcal{M}$ have no common vertex, i.e. whenever $u\!\leftarrow\!v,x\!\leftarrow\!y\!\in\!\mathcal{M}$ we have $|\{u,v,x,y\}|\!=\!4$;
	\item[$(2)$] for every edge $u\!\overset{w}{\leftarrow}\!v$ in $\mathcal{M}$, the corresponding weight $w$ is invertible in $R$;
	\item[$(3)$] denote by $\Gamma_{\!C_\ast}^{\mathcal{M}}$ the graph, obtained from $\Gamma_{\!C_\ast}$ by replacing every $u\!\overset{w}{\leftarrow}\!v\in\mathcal{M}$ with $u\!\overset{\!-\!1/\!w}{\longrightarrow}\!v$; then $\Gamma_{\!C_\ast}^\mathcal{M}$ must contain no directed cycles, and no infinite paths $u_1\!\to\!v_1\!\to\!u_2\!\to\!v_2\!\to\!\ldots$ with all $u_1,u_2,\ldots\!\in\!I_k$.
\end{enumerate}
\noindent Let $\mathcal{M}$ be a Morse matching. Let $\mathcal{M}_k\!=\!\mathcal{M}\!\cap\!\partial_k$ be the edges in $\mathcal{M}$ coming from $\partial_k$,
\begin{enumerate} \setlength{\itemindent}{-12mm}
	\item[] $I_k^+= \{v\!\in\!I_k;\, \exists u\!\leftarrow\!v\!\in\!\mathcal{M}\}= \{\text{indices of columns in }\partial_k\text{ that contain some }e\!\in\!\mathcal{M}\}$,
	\item[] $I_k^-= \{u\!\in\!I_k;\, \exists u\!\leftarrow\!v\!\in\!\mathcal{M}\}= \{\text{indices of rows in }\partial_{k+1}\text{ that contain some }e\!\in\!\mathcal{M}\}$,
	\item[] $I'_k= I_k\!\setminus\!(I_{k\!}^+\cup_{\!}I_k^-)\!=\! \{v\!\in\!I_k;\, v\text{ is not incident to any }e\!\in\!\mathcal{M}\}$, the \emph{critical} vertices.
\end{enumerate}
Let $\Gamma_{\!v,u}$ denote the set of all directed paths $\gamma$ in $\Gamma_{\!C_\ast}^{\mathcal{M}}$ from $v$ to $u$ (including paths of length 0 when $v\!=\!u$). Given such  $\gamma\!=(v\!\overset{w_0}{\to}\!u_1\!\!\overset{\!-\!1/\!w_1}{\longrightarrow}\!\! v_1 \!\overset{w_2}{\to}\!u_2\!\!\overset{\!-\!1/\!w_3}{\longrightarrow}\!\! v_2\!\to\!\ldots\! \longrightarrow\!\!v_r\overset{w_{2r}}{\to}\!\! u)$, let $\partial_\gamma(v)=w_{2r}\cdots\frac{-1}{w_3}w_2\frac{-1}{w_1}w_0\,u$ be the multiple of $u$ by the product of weights. Let $C'_k$ be the free module on $I'_k$. Define $\partial'_k\!\!: C'_k\!\to\!C'_{k_{\!}-\!1}$,\; $f_k\!\!: C'_k\!\to\!C_k$,\; $g_k\!\!: C_k\!\to\!C'_k$ by \[
\partial'_k(v')=\!\!\!\sum_{\substack{u'\in I'_{k\!-\!1},\\\gamma\in\Gamma_{\!v'\!,u'}}}\!\!\!\partial_\gamma(v'),\hspace{20pt}
f_k(v')=\!\!\!\sum_{\substack{v\in I_k,\\\gamma\in\Gamma_{\!v'\!\!,v}}}\!\!\!\partial_\gamma(v'),\hspace{20pt}
g_k(v)=\!\!\!\sum_{\substack{v'\in I'_k,\\\gamma\in\Gamma_{\!v,v'}}}\!\!\!\partial_\gamma(v). \tag{4} \label{(4)}\]
Let us recall the basic notions of homological algebra. A \emph{chain map} $\varphi_\ast\!:C_\ast\!\to\!D_\ast$ is a family of $R$-module morphisms $\varphi_k\!:C_k\!\to\!D_k$, such that $\varphi_{k-\!1}\partial_k^C \!=\! \partial_k^D\varphi_k$. In this case, $\varphi_\ast$ induces $R$-module morphisms $H_k\varphi\!:H_kC_\ast\!\to\!H_kD_\ast$, where $H_kC_\ast\!=\!\frac{\Ker\,\partial_k^C}{\Im\,\partial_{k+1}^C}$ is the \emph{$k$-th homology} of $C_\ast$ and $H_k\varphi$ sends $[x]\!\mapsto\![\varphi_k(x)]$. Chain maps $\varphi,\psi\!:C_\ast\!\to\!D_\ast$ are \emph{homotopic}, denoted $\varphi\!\simeq\!\psi$, if there exist morphisms $h_k\!:C_k\!\to\!D_{k+1}$ such that $\varphi_k\!=\!\psi_k+\partial^D_{k+1}h_k\!+\!h_{k-1}\partial^C_k$ for all $k$. In such a case, $H_\ast\varphi\!=\!H_\ast\psi$. A chain map $\varphi\!:C_\ast\!\!\to\!D_\ast$ is an \emph{h-equivalence} with \emph{h-inverse} $\psi\!:D_\ast\!\!\to\!C_\ast$, if $\varphi\!\circ\!\psi\simeq \id_{D_\ast}$ and $\psi\!\circ\!\varphi\simeq \id_{C_\ast}$. In such a case, $H_k\varphi\!=\!(H_k\psi)^{\!-\!1}\!\!:H_kC_\ast\!\to\!H_kD_\ast$ are isomorphisms. 
\par The fundamental result of AMT is stated as follows:

\begin{Thm}[\cite{articleSkoldbergMTFAV,articleJollenbeckADMTACA,articleKozlovDMTFCC} 2005]\label{AMT} \textbf{(a)} For any Morse matching $\mathcal{M}$, the induced $f_\ast\!:$ $(C'_{\!\ast},\partial'_{\!\ast}) \!\to\! (C_\ast,\partial_\ast)$ is an h-equivalence of chain complexes, with h-inverse $g_\ast$.\\
\textbf{(b)} There exist bases $\widetilde{I}_k$ for $C_k$, in which the boundary matrices have block forms\vspace{-1mm}
\[C_\ast\!:\ldots
\!\xleftarrow{\!\!\left[\!\begin{smallmatrix}\!\partial'_{_{\!}k\!-\!1}\!&0&0\\[-1pt]0&\!\!\id_{m_{\!k\!-\!1}}\!\!&0\\[-1pt]0&0&0\end{smallmatrix}\!\right]\!\!}
\!\!C_{\!k_{\!}-_{\!}1}\!\!\xleftarrow{\!\!\left[\!\begin{smallmatrix}\!\partial'_{_{\!}k}\!&0&0\\0&0&0\\0&0&\!\id_{_{\!}m_{\!k}}\!\end{smallmatrix}\!\right]\!\!}
\!\!C_{\!k}\!\!\xleftarrow{\!\!\left[\!\begin{smallmatrix}\!\partial'_{_{\!}k\!+\!1}\!&0&0\\[-1pt]0&\!\!\id_{m_{\!k\!+\!1}}\!\!&0\\[-2pt]0&0&0\end{smallmatrix}\!\right]\!\!}
\!\!C_{\!k+1}\!\!\xleftarrow{\!\!\left[\!\begin{smallmatrix}\!\partial'_{_{\!}k\!+\!2}\!&0&0\\[-1pt]0&0&0\\0&0&\!\id_{m_{\!k\!+\!2}}\!\end{smallmatrix}\!\right]\!\!} \!\ldots,\vspace{-1mm}\]
where $\id_r$ is the $r\!\times\!r$ identity matrix and $m_k\!=\!|\mathcal{M}_k|$.\\
\end{Thm}
\par Part (a) of \ref{AMT} holds even if $C_\ast$ is unbounded (i.e. $C_k\!\neq\!0$ for infinitely many $k\!\in\!\Z$) and $C_k$ is any (not necessarily finite) direct sum of submodules $\bigoplus_{v\in I_k}C_{k,v}$ over a (not necessarily commutative) unital ring. Then the column-finitary boundary matrices $\partial_k$ have for entries $\partial_{k,u,v}$ not weights $w\!\in\!R$ but morphisms $C_{k\!-\!1,u}\!\overset{\varphi}{\leftarrow}\!C_{k,v}$, condition $(2)$ says that $\varphi$ must be invertible as a morphism, and $\partial_\ast,f_\ast,g_\ast$ are defined via $\partial_\gamma= \varphi_{2r}\!\circ\!\ldots\!\circ\!(-\varphi_3^{\!-\!1}\!)\!\circ\! \varphi_2\!\circ\!(\!-\varphi_1^{\!-\!1}\!)\!\circ\!\varphi_0$. However, for our purposes, we shall work within the confines of assumptions stated in the conventions above. Hence the part of condition $(3)$ about no infinite paths is always satisfied.

\subsection{Example}\label{1.2.exp} Consider the Khovanov chain complex for the trefoil knot\vspace{-1mm} \[C_\ast\!:\: R^4\!
\xleftarrow{\left[\!\begin{smallmatrix}
	0&0&\gr{0}&\gr{0}&\gr{0}&\gr{0}\\
	1&0&\gr{\mns1}&\gr{0}&\gr{1}&\gr{0}\\
\rd{1}&0&\gr{\mns1}&\gr{0}&\gr{1}&\gr{0}\\
0&\rd{1}&\gr{0}&\gr{\mns1}&\gr{0}&\gr{1}\end{smallmatrix}\!\right]}\!R^6\!
\xleftarrow{\left[\!\begin{smallmatrix}
	\gr{\mns1}& 	       \gr{0}& 	  \gr{0}& \gr{0}&         \gr{1}& \gr{0}& \gr{0}& \gr{0}& \gr{0}& \gr{0}& \gr{0}& \gr{0}\\
	    \gr{0}&     \gr{\mns1}& \gr{\mns1}& \gr{0}& 	    \gr{0}& \gr{1}& \gr{1}& \gr{0}& \gr{0}& \gr{0}& \gr{0}& \gr{0}\\
\rd{\mns1}& 	   0& 	  0& 0& 	    0& 0& 0& \gr{0}& \gr{1}& \gr{0}& \gr{0}& \gr{0}\\
		0&\rd{\mns1}& \mns1& 0&         0& 0& 0& \gr{0}& \gr{0}& \gr{1}& \gr{1}& \gr{0}\\
		0& 	       0& 	  0& 0&\rd{\mns1}& 0& 0& \gr{0}& \gr{1}& \gr{0}& \gr{0}& \gr{0}\\
		0& 	       0& 	  0& 0&		    0&\rd{\mns1}&\mns1&  \gr{0}& \gr{0}& \gr{1}& \gr{1}& \gr{0}\end{smallmatrix}\!\right]}\!R^{12}\!
\xleftarrow{\left[\!\begin{smallmatrix}
	 \gr{1}&  	 \gr{0}& \gr{0}&  \gr{0}&  \gr{0}& \gr{0}& \gr{0}& \gr{0}\\
 	 \gr{0}&  	 \gr{0}& \gr{1}&  \gr{0}&  \gr{0}& \gr{0}& \gr{0}& \gr{0}\\
 	 0& 	 1& 0&  0& 	1& 0& 0& 0\\
 	 0&  	 0& 0&  1&  0& 0& 1& 0\\
	 \gr{1}& \gr{0}& \gr{0}&  \gr{0}&  \gr{0}& \gr{0}& \gr{0}& \gr{0}\\
 	 \gr{0}& \gr{1}& \gr{0}&  \gr{0}&  \gr{0}& \gr{0}& \gr{0}& \gr{0}\\
 	 0&  	 0& 1&  0& 	1& 0& 0& 0\\
 	 0&  	 0& 0&\rd{1}&0&1& 0& 0\\
\rd{1}&  	 0& 0&  0&  0& 0& 0& 0\\
 	 0& \rd{1}& 1&  0&  0& 0& 0& 0\\
 	 0&  	 0& 0&  0&\rd{1}& 0& 0& 0\\
 	 0&  	 0& 0&  0&0&\rd{1}& 1& 0
\end{smallmatrix}\!\right]}\!R^8.\vspace{-1mm}\]
Let the Morse matching consist of the red entries/edges. The associated graph $\Gamma_{C_\ast}$ (without displayed weights, and with $I_k\!=\!\{e_{k,1},e_{k,2},\ldots\}$) is below left. Gray entries / dotted edges are the ones that can be removed by remark \ref{1.3.rmk}\,(c). We omitted the arrows, since black and dotted lines always point leftwards, red lines point rightwards. Below right are the paths between critical vertices, that give $\partial'_\ast$.\vspace{-2mm} \[\hspace{-2mm}
\begin{tikzpicture}\matrix (m) [matrix of math nodes, row sep=-1pt, column sep=25pt]{
\ccircled{$e_{0,1}$}&    e_{1,1}& e_{2,1}				& e_{3,1}\\
					&  			& e_{2,2}				& e_{3,2}\\
					& 	 e_{1,2}& \ccircled{$e_{2,3}$}	& \\
\ccircled{$e_{0,2}$}&  			& \ccircled{$e_{2,4}$}	& \ccircled{$e_{3,3}$}\\
					& 	 e_{1,3}& e_{2,5}				& e_{3,4}\\
					&  			& e_{2,6}				& \\
e_{0,3}				& 	 e_{1,4}& \ccircled{$e_{2,7}$}	& e_{3,5}\\
					&  			& e_{2,8}				& e_{3,6}\\
					&    e_{1,5}& e_{2,9}				& \\
e_{0,4}				&  			&e_{2,10}				& \ccircled{$e_{3,7}$}\\
					&    e_{1,6}&e_{2,11}				& \ccircled{$e_{3,8}$}\\
					&  			&e_{2,12}				&  \\};
\draw(m-4-1.east)--(m-1-2.west); 	  \draw[dotted](m-4-1.east)--(m-5-2.west); \draw[dotted](m-4-1.east)--(m-9-2.west);
\draw[red](m-7-1.east)--(m-1-2.west); \draw[dotted](m-7-1.east)--(m-5-2.west); \draw[dotted](m-7-1.east)--(m-9-2.west);
\draw[red](m-10-1.east)--(m-3-2.west);\draw[dotted](m-10-1.east)--(m-7-2.west);\draw[dotted](m-10-1.east)--(m-11-2.west);
\draw[dotted](m-1-2.east)--(m-1-3.west); \draw[dotted](m-1-2.east)--(m-5-3.west);
\draw[dotted](m-3-2.east)--(m-2-3.west); \draw[dotted](m-3-2.east)--(m-3-3.west); \draw[dotted](m-3-2.east)--(m-6-3.west); \draw[dotted](m-3-2.east)--(m-7-3.west);
\draw[red](m-5-2.east)--(m-1-3.west); \draw[dotted](m-5-2.east)--(m-9-3.west);
\draw[red](m-7-2.east)--(m-2-3.west); \draw(m-7-2.east)--(m-3-3.west); \draw[dotted](m-7-2.east)--(m-10-3.west); \draw[dotted](m-7-2.east)--(m-11-3.west);
\draw[red](m-9-2.east)--(m-5-3.west); \draw[dotted](m-9-2.east)--(m-9-3.west);
\draw[red](m-11-2.east)--(m-6-3.west); \draw(m-11-2.east)--(m-7-3.west); \draw[dotted](m-11-2.east)--(m-10-3.west); \draw[dotted](m-11-2.east)--(m-11-3.west);
\draw[dotted](m-1-3.east)--(m-1-4.west);
\draw[dotted](m-2-3.east)--(m-4-4.west);
\draw(m-3-3.east)--(m-2-4.west); \draw(m-3-3.east)--(m-7-4.west);
\draw(m-4-3.east)--(m-5-4.west); \draw(m-4-3.east)--(m-10-4.west);
\draw[dotted](m-5-3.east)--(m-1-4.west);
\draw[dotted](m-6-3.east)--(m-2-4.west);
\draw(m-7-3.east)--(m-4-4.west); \draw(m-7-3.east)--(m-7-4.west);
\draw[red](m-8-3.east)--(m-5-4.west); \draw(m-8-3.east)--(m-8-4.west);
\draw[red](m-9-3.east)--(m-1-4.west);
\draw[red](m-10-3.east)--(m-2-4.west); \draw(m-10-3.east)--(m-4-4.west);
\draw[red](m-11-3.east)--(m-7-4.west);
\draw[red](m-12-3.east)--(m-8-4.west); \draw(m-12-3.east)--(m-10-4.west);
\end{tikzpicture}                       \hspace{12mm}
\begin{tikzpicture}\matrix (m) [matrix of math nodes, row sep=-3pt, column sep=15pt]{ \gamma\!:\\
	e_{2,10}			&\ccircled{$e_{3,3}$}\\
	\ccircled{$e_{2,3}$}&e_{3,2}\\[5mm]
	\ccircled{$e_{2,7}$}&\ccircled{$e_{3,3}$}\\[5mm]
	\ccircled{$e_{2,4}$}&\ccircled{$e_{3,7}$}\\[5mm]
	e_{2,12}			&\ccircled{$e_{3,7}$}\\
	e_{2,8}				&e_{3,6}\\
	\ccircled{$e_{2,4}$}&e_{3,4}\\};
\draw(m-2-1.east)--(m-2-2.west); \draw[red](m-2-1.east)--(m-3-2.west); \draw(m-3-1.east)--(m-3-2.west);
\draw(m-4-1.east)--(m-4-2.west); \draw(m-5-1.east)--(m-5-2.west);
\draw(m-6-1.east)--(m-6-2.west); \draw[red](m-6-1.east)--(m-7-2.west); \draw(m-7-1.east)--(m-7-2.west);
\draw[red](m-7-1.east)--(m-8-2.west); \draw(m-8-1.east)--(m-8-2.west);
\end{tikzpicture}
\begin{tikzpicture}\matrix (m) [matrix of math nodes, row sep=-3pt, column sep=15pt]{ \partial_\gamma\!:\\[3mm]
	1\!\cdot\!\frac{-\!1}{1}\!\cdot\!1=-1\\[6mm]
	1\\[7mm]
	1\\[14.5mm]
	1\!\cdot\!\frac{-\!1}{1}\!\cdot\!1\!\cdot\!\frac{-\!1}{1}\!\cdot\!1=1\\	};
\end{tikzpicture}\vspace{-2mm}\]
From this, we get an h-equivalent complex $C'_\ast\!: R^2\!
\xleftarrow{\,0\,}\!R^0\!
\xleftarrow{\,0\,}\!R^3\!
\xleftarrow{\!\!\left[\!\begin{smallmatrix}\mns1\!&0&0\\0&2&0\\1&0&0\end{smallmatrix}\!\right]\!\!}\!R^3\!$.
Moreover, $f_\ast\!=\!\left[\!\begin{smallmatrix}1&0\\0&1\\0&0\\0&0\end{smallmatrix}\!\right]\!\!, 0,\! \left[\!\begin{smallmatrix}0&0&0\\\mns1&0&0\\1&0&0\\0&1&0\\0&0&0\\0&0&\mns1\\0&0&1\\0&0&0\\0&0&0\\0&0&0\\0&0&0\\0&0&0\end{smallmatrix}\!\right]\!\!,\! \left[\!\begin{smallmatrix}0&0&0\\\mns1&0&0\\1&0&0\\0&1&0\\0&0&0\\0&\mns1&0\\0&1&0\\0&0&1\\\end{smallmatrix}\!\right]$ and $g_\ast\!=\!
\left[\!\begin{smallmatrix}1&0&0&0\\0&1&\mns1&0\end{smallmatrix}\!\right]\!\!, 0,\!
\left[\!\begin{smallmatrix}0&0&1&0&0&0&0&0&0&\mns1&\mns1&0\\
0&0&0&1&0&0&0&\mns1&0&0&0&1\\ 0&0&0&0&0&0&1&0&0&0&\mns1&0\end{smallmatrix}\!\right]\!\!,
\left[\!\begin{smallmatrix}0&0&1&0&0&0&0&0\\0&0&0&0&0&0&1&0\\0&0&0&0&0&0&0&1\end{smallmatrix}\!\right]\!$. Thus over $R\!=\!\Z$, the homology and its generators are $H_0C_\ast\!=\! \langle e_{0,1},e_{0,2}\rangle \!\cong\!\Z^2\!$, $H_1C_\ast\!\cong\!0$, $H_2C_\ast\!=\! \langle e_{2,3}\!-\!e_{2,2},e_{2,4}\rangle \!\cong\!\Z\!\oplus\!\Z_2$, $H_3C_\ast\!=\!\langle e_{3,8}\rangle\!\cong\!\Z$. \hfill$\lozenge$\vspace{2mm}

\subsection{Remarks}\label{1.3.rmk} The following observations are crucial for our algorithm \ref{5.1.algorithm}.\\[3pt]
\textbf{\underline{\smash{(a)}}} If $C_\ast$ is the Poincar\'{e} chain complex of a finite simplicial complex $\Delta$ (or regular CW-complex), then $\Gamma_{\!C_\ast}$ is the Hasse diagram of $\Delta$ and AMT reduces to Forman's DMT \cite{articleFormanMTCC} (discrete Morse theory, 1995). However, DMT produces a topological h-equivalence of CW-complexes $\Delta\!\simeq\!\Delta'$, not just an algebraic h-equivalence $C_\ast\!\simeq\!C'_\ast$. To determine the gluing maps, one needs to take into account paths between all critical simplices, not just those in consecutive degree. See \cite{articleNandaDMTL} for an elaboration. \\[5pt]
\textbf{\underline{\smash{(b)}}} Notice that $\partial'_k$ is determined by $\partial_{k-_{\!}1},\partial_k,\partial_{k+1}$, not just $\partial_k$. Indeed, $I'_k$ is obtained from $\partial_{k-_{\!}1},\partial_k$ and $I'_{k+1}$ is obtained from $\partial_k,\partial_{k+1}$, whilst $\partial'_k$ is computed from $I'_k,I'_{k+1},\mathcal{M}_k,\partial_k$. Also, (\ref{(4)}) has a more precise formulation, namely \vspace{-1mm}
\[\textstyle{f_k(v')=v'\!+\! \sum_{v\in I_k^+\!,\gamma\in\Gamma_{\!v'\!\!,v}}\!\partial_\gamma(v')} \text{ \;\;\;and\;\;\; } g_k(v)=
\Bigg\{\begin{smallmatrix}
v &\text{\!\!; if }v\in I'_k\\
0 &\text{\!\!; if }v\in I^+_k\\
\sum_{v'\in I'_k,\gamma\in\Gamma_{\!v,v'}}\!\partial_\gamma(v) &\text{\!\!; if }v\in I^-_k.\end{smallmatrix}\vspace{-1mm}\]
Hence, $\partial'_k$ and $f_k$ contain only edges from $\partial_k$, whilst $g_k$ has only edges from $\partial_{k+1}$. Also, $f_0$ is an inclusion ($I'_0$-columns of $\id_{|I_0|}$) and $g_N$ a projection ($I'_N$-rows of $\id_{|I_N|}$).\\[5pt]
\textbf{\underline{\smash{(c)}}} Let us see how AMT allows many matrix entries to be deleted and some to be edited before computation even begins, thereby lowering space and time complexity. 
\par For every $u\!\leftarrow\!v$ in $\mathcal{M}_k$, no edge $x\!\leftarrow\!u$ or $v\!\leftarrow\!y$ in $\Gamma_{C_\ast}^{\mathcal{M}}$ can lie on any path $\gamma$ from $\partial'_\ast$ or $f_\ast$ or $g_\ast$, since $\mathcal{M}$ is a matching. Thus we can set to zero all entries in the $u$-column (or $I_{k-\!1}^-$-columns) of $\partial_{k-\!1}$ and $v$-row (or $I_k^+$-rows) of $\partial_{k+1}$, without changing the computed complex $C'_\ast$ and h-equivalences $f_\ast,g_\ast$ (using the same $\mathcal{M}$). This modification has a useful consequence: when calculating $\partial'_k$, any partially computed zig-zag path either continues (its endpoint is in $I_{k-1}^-$) or finishes (its endpoint is in $I'_{k-1}$). In other words, there are no dead-end paths; not having to check every time whether the endpoint is in $I_{k-1}^+$ makes digraph traversal faster. 
\par On a related note, for every $u\!\overset{w}{\leftarrow}\!v$ in $\mathcal{M}_k$, we can replace the basis element $u$ with $-w\,u$, i.e. multiply the $u$-row of $\partial_k$ by $\frac{-\!1}{w}$ (and $u$-column of  $\partial_{k-1}$ by $-w$, but that is now zero, by the above). Then all edges of $\mathcal{M}$ in $\Gamma_{C_\ast}^{\mathcal{M}}$ have weight $1$, which simplifies the formula above (\ref{(4)}) to $\partial_\gamma(v)=w_{2r}\!\cdots w_4w_2w_0\,u$. This way, we avoid computing the same inverse weight several times in every call of an iterative/recursive function for finding paths. From the altered $\widetilde{\partial}_k$, we compute $\widetilde{\partial}'_k,\widetilde{f}_k,\widetilde{g}_{k-1}$, but by the nature of paths in $\Gamma_{C_\ast}^{\mathcal{M}}$, we have $\widetilde{\partial}'_k\!=\!\partial'_k$, $\widetilde{f}_k\!=\!f_k$, and $\widetilde{g}_{k-1}$ equals $g_{k-1}$ with each $u$-column multiplied by $-w$. Thus, we multiply each $u$-column of $\widetilde{g}_{k-1}$ by $\frac{-\!1}{w}$ and obtain $\partial'_k,f_k,g_{k-1}$ w.r.t. original bases.\\[5pt]
\textbf{\underline{\smash{(d)}}} Let us see how AMT produces generators in homology. As long as one of the matrices $\partial_k$ contains at least one invertible entry $\partial_{k,u,v}\!=\!w$, we always have a nonempty Morse matching. Indeed, just take $\mathcal{M}\!=\!\{u\!\overset{w}{\leftarrow}\!v\}$. But for more efficiency, picking $\mathcal{M}_\leq$ from \ref{2.1.formulation} relative to a good order $\leq$ according to \ref{3.ordering}, is advisable. Hence there is a sequence of h-equivalences\vspace{-2mm}
\[C_\ast \stackrel[g'_\ast]{f'_\ast}\longleftrightarrows C'_\ast \stackrel[g''_\ast]{f''_\ast}\longleftrightarrows C''_\ast \longleftrightarrows\ldots \stackrel[g^{(r)}_\ast]{f^{(r)}_\ast}\longleftrightarrows  C^{(r)}_\ast\!\!, \tag{5} \label{(5)}\vspace{-2mm}\]
where all entries of $\partial_\ast^{(r)}$ are nonunits. If $R$ is a field, we have $\partial_\ast^{(r)}\!\!=\!0$, so $C^{(r)}_k\!\cong\!H_k C_\ast$; then the image of $f_\ast\!:=\! f_\ast'\!\circ\!f_\ast''\!\circ\!\ldots\!\circ\!f_\ast^{(r)}$ gives the homology generators inside $C_\ast$. The h-inverse of $f_\ast$ is then $g_\ast\!:=\!g^{(r)}_\ast\!\circ\!\ldots\!\circ\!g''_\ast\!\circ\!g'_\ast$, i.e. $g_\ast\!\circ\!f_\ast\!=\!\id_{C_\ast^{(r)}}$ and $f_\ast\!\circ\!g_\ast\!\simeq\!\id_{C_\ast}$. A special case of this result was presented in \cite{articleKaczynskiMrozekSlusarekHCRCC} (when $\mathcal{M}$ consists of only one edge) and \cite{articleDlotkoWagnerCHPHIMD} (when $\mathcal{M}$ is a specific algorithmically constructed matching). A similar but non-sparse approach using contractions was employed in \cite{articlePilarczykRealCCHCCOCC}. \\[5pt]
\textbf{\underline{\smash{(e)}}} Let us see how AMT reveals torsion in homology. If $R_{(p)}$ is a localization of a PID $R$ at a prime $p$, then the noninvertible entries are multiples of $p$. Thus $\smash{\partial_\ast^{(r)}}\!=\!p^{a}\widetilde{\partial}_\ast^{(r)}$ for some $a\!\in\!\N$ and $\widetilde{\partial}_\ast^{(r)}$. Via \ref{AMT}\,(b) we reconstruct $H_\ast(C^{(r)}_\ast\!\!,\partial^{(r)}_\ast)$ from $H_\ast(C^{(r)}_\ast\!\!,\smash{\frac{\partial^{(r)}_\ast}{p^a}})$. Hence we continue with (\ref{(5)}) on \smash{$(C^{(r)}_\ast\!\!,\frac{\partial^{(r)}_\ast}{p^a})$}: for every Morse matching $\mathcal{M}$ on it, we add \smash{$(\frac{R}{Rp^{a}})^{|\mathcal{M}_k|}$} to $H_{k_{\!}-\!1}(C_\ast,\partial_\ast)$. This process ends at $\partial_\ast^{(r_{1}\!+\ldots+r_{t})}\!\!=\!0$, and we obtain the $p$-torsion of $H_{\!\ast}(C_{\!\ast};_{\!}R)$. See \ref{4.nonunits} for a more detailed account. \\[5pt]
\textbf{\underline{\smash{(f)}}} Let us see how AMT induces morphisms on homologies. Let $\varphi\!:C_\ast\!\to\!D_\ast$ be a chain map and $R$ a field. Then (\ref{(5)}) produces chain complexes $C'_\ast$ and $D'_\ast$ with zero boundaries. Since $H_k(C_\ast)\!\cong\!C'_k$ and $H_k(D_\ast)\!\cong\!D'_k$ for all $k$, the h-equivalences reveal the induced maps in the form of composition \smash{$H_\ast(\varphi)\!=\!\varphi': C'_\ast\!\overset{f_C}{\to}\! C_\ast\!\overset{\varphi}{\to}\! D_\ast\!\overset{g_D}{\to}\!D'_\ast$}.\\[5pt] 
\textbf{\underline{\smash{(g)}}} Let us see how AMT can compute relative homology. Suppose $A_\ast\leq B_\ast$ is a chain subcomplex over a field. Let $I_k$ (resp. $J_k$) be a basis of $A_k$ (resp. $B_k$) and suppose that $I_k\!\subseteq\!J_k$, for all $k$. This condition is often satisfied, e.g. for any CW subcomplexes (Poincar\'{e} complex), any subgroups (Eilenberg-MacLane complex), as well as many associative subalgebras (Hochschild complex) and Lie subalgebras (Chevalley complex), we indeed have $I_k\!\subseteq\!J_k$. To compute $H_\ast(B_\ast/A_\ast)$, we delete all $I_{k-1}$-rows and $I_k$-columns from $\partial^B_k$, for all $k$, and continue with AMT on the smaller complex. However, if $I_k\!\nsubseteq\!J_k$, then we need to find sparse bases for $B_k/A_k$ and express $\partial^B_k$ in those bases. This is achieved by AMT using the results in \ref{5.4.SparseLinearAlgebra}.\\[5pt]
\textbf{\underline{\smash{(h)}}} Let us see how AMT provides (co)homology operations. Since a cochain complex is formally the same as a chain complex, only the grading goes in the other direction, AMT can be applied. Over a field, any finite-dimensional complexes $C_\ast$ and $C^\ast$ are reduced to their (co)homology using (\ref{(5)}). Then any operations $C_i\!\otimes\!C_j\!\to\!C_k$, $C^i\!\otimes\!C^j\!\to\!C^k$, $C_i\!\otimes\!C^j\!\to\!C_k$, $C^i\!\otimes\!C_j\!\to\!C^k$, or similar, that descend onto (co)homology, can be read-off using the method in (f). A non-sparse approach to this was employed in \cite{articleDiazRealCCOFSC} using contractions via SNFs. \\[5pt]
\textbf{\underline{\smash{(i)}}} Let us see how DMT can only do the first step of AMT. In general, it is necessary to use AMT several times on a complex. More precisely, for most $C_\ast$ there does not exist a single Morse matching that would reduce it to a complex with zero boundaries. For instance, we have two reductions \vspace{-2mm} \[
C_\ast\!: \R^2\!\!\xleftarrow{\left[\begin{smallmatrix}1&1&0\\\rd{1}&-\!1&0\end{smallmatrix}\right]}\!\!\R^3\hspace{15pt} \xleftarrow{f^{(1)}}\hspace{15pt} 
C^{(1)}_\ast\!\!: \R^1\!\!\xleftarrow{\left[\begin{smallmatrix}\rd{2}&0\end{smallmatrix}\right]}\!\!\R^2\hspace{15pt} \xleftarrow{f^{(2)}}\hspace{15pt} 
C^{(2)}_\ast\!\!: \R^0\!\xleftarrow{0}\R^1,\vspace{-2mm}\] but this cannot be achieved in just one step, since a Morse matching on this $C_\ast$ can have at most one element. Even if $C_\ast$ comes from a simplicial complex, in most cases a perfect Morse function does not exist. Hence AMT has a big advantage over DMT, in that the latter can perform only the first reduction of the former. \vspace{4mm}

\begin{Lmm}\label{lmm0} If we apply AMT (relative to any nonempty Morse matching) repeatedly until the boundary matrices are zero, then the resulting rule $C_\ast\!\mapsto\!C'_\ast$ and $\varphi\!\mapsto\!\varphi'$ from \ref{1.3.rmk}(d,f), is the homology functor from the category of finite-dimensional chain complexes of vector spaces to the category of $\Z$-graded vector spaces.
\end{Lmm}\vspace{-1mm}
Hence any commuting diagram of chain maps is transformed (with repeated applications of AMT) into a commuting diagram of homology maps.
\begin{proof} We wish to show that $\id_{C_\ast}'=\id_{C'_\ast}$ and $(\psi_\ast\varphi)'=\psi_\ast'\varphi_\ast'$. Notice that in general, if complexes $A_\ast$ and $B_\ast$ have zero boundary maps, then any chain maps $\varphi_\ast,\psi_\ast\!:A_\ast\!\to\!B_\ast$ are homotopic iff they are equal. 
\par Suppose that a sequence of Morse matchings on $C_\ast$ results in an h-equivalence: \vspace{-2mm} \[\xymatrix{ 
	C_\ast \ar[r]^{\id_\ast} \ar@<2pt>[d]^{g_\ast}	& C_\ast \ar@<2pt>[d]^{g_\ast}\\
	C'_\ast \ar@<2pt>[u]^{f_\ast}					& C'_\ast \ar@<2pt>[u]^{f_\ast} }\vspace{-1mm}\]
\noindent Then $g_\ast f_\ast\!\simeq\!\id_{C'_\ast}$ and $f_\ast g_\ast\!\simeq\!\id_{C_\ast}$\!. Since $C'_\ast$ has zero boundaries, $\id_{C_\ast}'\!\!=\!g_\ast f_\ast\!\!=\!\id_{C'_\ast}$.
\par Suppose that three sequences of Morse matchings result in h-equivalences: \vspace{-2mm} \[\xymatrix{ 
	C_\ast \ar[r]^{\varphi_\ast} \ar@<2pt>[d]^{g^C_\ast}& D_\ast \ar@<2pt>[d]^{g^D_\ast} \ar[r]^{\psi_\ast} & E_\ast \ar@<2pt>[d]^{g^E_\ast}\\
	C'_\ast \ar@<2pt>[u]^{f^C_\ast}						& D'_\ast \ar@<2pt>[u]^{f^D_\ast}					& E'_\ast \ar@<2pt>[u]^{f^E_\ast} }\vspace{-1mm}\] 
We wish to show that $g^E_\ast\psi_\ast\varphi_\ast f^C_\ast=g^E_\ast\psi_\ast f^D_\ast g^D_\ast\varphi_\ast f^C_\ast$. We have $f^D_\ast g^D_\ast\!\simeq\!\id_{D_\ast}$ via a homotopy $h^D_\ast\!: D_{\ast-\!1}\!\!\to\!D_\ast$, meaning that $f^D_kg^D_k\!=\!\id_{D_k}\!\!+\!h^D_k\partial_k^D\!+\!\partial_{k+1}^Dh_{k+1}^D$ for all $k$. Thus it suffices to show that 
$g^E_\ast\psi_\ast h^D_\ast\partial_\ast^D\varphi_\ast f^C_\ast =0= g^E_\ast\psi_\ast \partial_{\ast+1}^Dh_{\ast+1}^D\varphi_\ast f^C_\ast$. To this end, $\partial_\ast^D\varphi_\ast f^C_\ast =\varphi_{\ast-\!1}\partial_\ast^C f^C_\ast =\varphi_{\ast-\!1}f^C_{\ast-\!1} {\partial'_\ast}^{\!C} =0$ and $g^E_\ast\psi_\ast \partial_{\ast+1}^D= g^E_\ast\partial_{\ast+1}^E\psi_{\ast+1}= {\partial'^E}_{\!\!\!\!\!\ast+1} g^E_{\ast+1}\psi_{\ast+1}= 0$, since $f_\ast^C\!,g^E_\ast\!,\varphi_\ast,\psi_\ast$ are chain maps and ${\partial'_\ast}^{\!C}\!\!=\!0\!=\!{\partial'_\ast}^{\!E}$.
\end{proof}\vspace{4mm}

\section{Matchings from orderings}
In this section, we define a class of matchings on \emph{any} chain complex of free finite-rank $R$-modules. The construction is dependent on the choice of ordering of basis elements. This makes AMT applicable to a vast array of (co)homology theories. 

\subsection{Formulation}\label{2.1.formulation} Pick a total order $\leq_k$ on $I_k$ and let $\leq\,=\bigcup_k\!\leq_k$ be the combined partial order on the vertices of $\Gamma_{C_\ast}$. Visualize the elements of $I_k$ positioned vertically, with $u$ above $v$ iff $u\!<\!v$ (as with row indices in the matrix $\partial_{k+1}$). Define
\[\mathcal{M}_\leq=\Big\{u\!\overset{w}{\leftarrow}\!v;\, \begin{smallmatrix}\forall x\leftarrow v:\: x<u,\\ \forall u\leftarrow y:\: v<y\end{smallmatrix}\Big\}\cap
\Big\{u\!\overset{w}{\leftarrow}\!v;\, w\!\in\!R^{\times}\Big\},\] the set of all steepest edges which also have invertible weights. The critical vertices are $I'_\ast\!=\!\Big\{v;
\begin{smallmatrix}\text{if }v\overset{r}{\leftarrow} u\text{ is steepest into }v\text{,\;\;\;\; then }r\in R\setminus R^\times\!\text{ or }\exists v'\!\leftarrow u\text{ with }v'>v,\\
                   \text{if }w\overset{s}{\leftarrow} v\text{ is steepest out of }v\text{, then }s\in R\setminus R^\times\!\text{ or }\exists w\leftarrow v'\text{ with }v'\!<v\end{smallmatrix}\Big\}$.
Visually, \[\mathcal{M}_\leq\!=\!\Big\{ \raisebox{-12pt}{\includegraphics[width=0.17\textwidth]{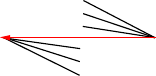}}\Big\}\text{ \;and\;
}I'_\ast\!=\! \Big\{v; \raisebox{-10pt}{\includegraphics[width=0.22\textwidth]{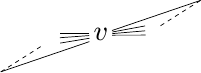}}\Big\}.\]

\begin{Lmm}\label{lmm1} \!\textbf{(a)}\! $\mathcal{M}_\leq$ is a Morse matching, the \emph{steepness} pairing associated to $\leq$.\\
\textbf{(b)} Every Morse matching is a subset of some steepness matching, so \vspace{-3pt}
\[\{\text{maximal Morse matchings on }\Gamma_{C_\ast}\} \subseteq \{\text{steepness pairings on }\Gamma_{C_\ast}\}.\] 
\end{Lmm}\vspace{-1mm}
\begin{proof} \textbf{(a)} The three conditions from the above definition must be verified.
\par $(1)$ $\mathcal{M}$ is a matching: For any edges
\!\!\raisebox{4pt}{$\xymatrixrowsep{-7pt}\xymatrixcolsep{14pt}
\xymatrix{&\scriptstyle x\ar[ld]\\\scriptstyle u&\\ &\scriptstyle y\ar[lu]\\}$}\!\!, either $x\!<\!y$ or $y\!<\!x$, so both edges cannot be in $\mathcal{M}$. For any edges \!\!\raisebox{4pt}{\smash{$\xymatrixrowsep{-7pt}\xymatrixcolsep{14pt}
\xymatrix{\scriptstyle x&\\&\scriptstyle v\ar[ld]\ar[lu]\\ \scriptstyle y&\\}$}}\!\!, either $x\!<\!y$ or $y\!<\!x$, so both edges cannot be in $\mathcal{M}$. For any edges $x\!\overset{r}{\leftarrow}\!y\!\overset{s}{\leftarrow}\!z$, either $rs\!=\!0$ (then weights $r$ and $s$ are not units of $R$, so their edges cannot be in $\mathcal{M}$), or there exist edges $x\!\leftarrow\!v\!\leftarrow\!z$ with $y\!\neq\!v$ (otherwise $\partial_{k\!-\!1}\partial_k\!\neq\!0$). Now, either $y\!<\!v$ or $v\!<\!y$, so $x\!\leftarrow\!y$ and $y\!\leftarrow\!z$ cannot be both in $\mathcal{M}$, because at least one of them would not be the steepest.
\par $(2)$ $\mathcal{M}$ consists of isomorphisms: This is by the assumption on invertible weights.
\par $(3)$ $\mathcal{M}$ induces no cycles: Every zig-zag \!\!\raisebox{5pt}{$\xymatrixrowsep{-6pt}\xymatrixcolsep{12pt}
\xymatrix{&\scriptstyle x\ar[ld]\\ \scriptstyle u&\\ &\scriptstyle y\ar[lu]\\}$}\!\! in which $u\!\leftarrow\!x\in\mathcal{M}$ implies that $x\!<\!y$. Thus every zig-zag strictly increases the initial vertex, so a zig-zag path cannot end at its starting point, because it would imply $y\!>\!y$.\\
\textbf{(b)} Let $\mathcal{M}$ be an arbitrary Morse matching on $\Gamma_{C_\ast}$. It suffices to find for every $k$ a total order $\leq_k$ on $I_k$ such that the edges in $\mathcal{M}$ are the steepest.
\par By assumption, the digraph $\Gamma^{\mathcal{M}}_{C_\ast}$ is acyclic, so it can be viewed as a poset: $u\!\preceq\!v$ iff there exists a directed path from $u$ to $v$. Let $\preceq_k$ be the restriction of $\preceq$ to $I_k$. By construction, members of $\mathcal{M}$ are steepest w.r.t. $\preceq_k$. Since every finite partial order $\preceq_k$ can be extended to a linear order $\leq_k$, we conclude that $\mathcal{M}\!\subseteq\!\mathcal{M}_\leq$.
\par If $\mathcal{M}$ is also maximal, then $\mathcal{M}\!\subseteq\!\mathcal{M}_\leq$ for some $\leq$ implies $\mathcal{M}\!=\!\mathcal{M}_\leq$.
\end{proof}
Note that if $\partial_k$ was infinite, with only diagonal and first supdiagonal entries which were all units, then $\mathcal{M}_\leq$ would not satisfy the part of condition $(3)$ about infinite paths. Thus infinite steepness pairings are not always Morse matchings.

\subsection{Matrix interpretation} If $\partial_k$ are given as finite matrices, we may assume that $I_k\!=\!\{1,\ldots,|I_k|\}$ and $\leq_k$ is the usual total order on $\N$, so the first row/column index is the smallest and the last one is the largest. Our $\mathcal{M}_\leq$ consists of those invertible matrix entries that have only zeros left in its row and below in its column, \vspace{-2mm}
\[\mathcal{M}_\leq=\left\{\raisebox{4.5pt}{$\scriptstyle u\,$}\!\!
\overset{\hspace{4pt}v}{\left[\hspace{-3pt}\begin{smallmatrix}~\\~\\~\\0\cdots0\rd{w}~~~\\
\hspace{15pt}0~~~\\[-1pt]
\hspace{14.5pt}\scalebox{0.8}{\rotatebox{90}{$.\hspace{0.4pt}.\hspace{0.4pt}.$}}~~~\\[-0.5pt]
\hspace{15pt}0~~~\\[2pt]\end{smallmatrix}\hspace{3pt}\right]}\!;\, w\!\in\!R^\times\!\!\right\}\!.
\tag{6} \label{(6)}\vspace{-1mm}\]
Thus $\mathcal{M}$ is largest when the matrices $\partial_k$ have 'block upper-triangular' form. The set of critical vertices is $I_k'=I_k\!\setminus(_{\!}\{\text{column indices of }\mathcal{M}_k\!\}\!\cup\!\{\text{row indices of }\mathcal{M}_{k+_{\!}1}\!\}_{\!})$.\\

\par Let us recall Gaussian elimination. Given a chain complex $C_\ast$ and a pivot (which is an invertible entry \smash{$u\!\overset{w}{\leftarrow}\!v\in\partial_k$}), for every nonzero entry \smash{$u'\!\overset{w'}{\leftarrow}\!v\in\partial_k$} (resp. \smash{$u\!\overset{w'}{\leftarrow}\!v'\in\partial_k$}) we add $\tfrac{-w'}{w}\!\cdot\!u$-row to the $u'$-row (resp. $\tfrac{-w'}{w}\!\cdot\!v$-column to the $v'$-column) in $\partial_k$. Since this operation alters the basis of $C_{k-1}$ (resp. $C_k$), namely replaces $u$ by $u\!+\!\tfrac{w'}{w}\!\cdot\!u'$ (resp. $v'$ by $v'\!-\!\tfrac{w'}{w}\!\cdot\!v$), we also have to add $\tfrac{w'}{w}\!\cdot\!u'$-column to the $u$-column in $\partial_{k-1}$ (resp. $\tfrac{w'}{w}\!\cdot\!v'$-row to the $v$-row in $\partial_{k+1}$). After this is done, we have changed all the entries in the $u$-row and $v$-column of $\partial_k$ (except the pivot) to zero. Unfortunately, we have also altered entries in other rows/columns of $\partial_{k-1},\partial_k,\partial_{k+1}$; this typically leads to an increase in matrix density, a major obstacle in sparse matrix computations. Nevertheless, since $u\!\overset{w}{\leftarrow}\!v$ is the only nonzero entry in its row and column, it represents a subcomplex in $C_\ast$ which is a direct summand, and since $w$ is invertible, this subcomplex is contractible. Consequently, we can remove the basis elements $u\!\in\!I_{k-1}$ and $v\!\in\!I_k$ (i.e. delete the $u$-row of $\partial_k$, $u$-column of $\partial_{k-1}$, $v$-column of $\partial_k$, $v$-row of $\partial_{k+1}$) and get a homotopy-equivalent chain complex. Repeating this over a field until there are no more unit entries leads to all matrices $\partial_k$ being zero, hence the new (smaller) $C_k$ is the homology. This algorithm is the content of the paper \cite{articleKaczynskiMrozekSlusarekHCRCC}. Of course, every row/column operation on a sparse matrix is costly, since it reshapes the whole data structure, so it is not advisable to be doing this pivot after pivot, but rather in tandem. Next, we show that AMT is a generalization of this procedure.

\begin{Lmm}\label{lmm2} Given $\mathcal{M}_\leq$, the $u'\!$-row  $r'$ of $\partial'_k$ is obtained from the $u'\!$-row $r$ of $\partial_k$ in the following way. While $V\!:=\!\{\text{nonzero positions of }r\}\!\cap\!I_k^+$ is nonempty, for $v\!=\!\min V$\!, $u\!\overset{w}{\leftarrow}\!v\in\mathcal{M}$, \smash{$u'\!\!\overset{\:w'}{\leftarrow}\!v\in\Gamma_{\!C_\ast}$}\!, add $\frac{-w'}{w}\!\cdot\!(u$-row of $\partial_k)$ to $r$. At the end, remove all $I_k^+$-entries and $I_{k+1}^-$-entries from $r$ to get the shorter row $r'$\!.
\end{Lmm}\vspace{-1mm}
\begin{proof} For the sake of clarity, we assume $\partial_k$ have block upper-triangular form, though the arguments work for any $\leq$. This is depicted by the image below: white area are zeros, full (vertical) lines are nonzeros, dashed (horizontal) lines and gray area are zeros and nonzeros, $\mathcal{M}_\leq$ consists of those red entries that are invertible, black bullets are nonzero entries, the dotted line is a path $\gamma$ from $v'$ to $u'$ in $\partial_k$.\vspace{-2mm}
\[\includegraphics[width=0.5\textwidth]{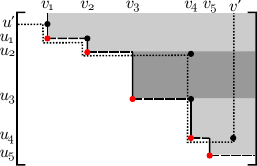} \tag{7} \label{(7)}\vspace{-1mm}\]
By the construction of $\Gamma_{\!C_\ast}^{\mathcal{M}}$, a path $\gamma\!\in\!\Gamma_{v'\!,u'\!}$ in $\partial_k$ consists of horizontal moves from black to red bullets, and vertical moves from red to black bullets. Let $I_k^+$ consist of $v_1,v_2,\ldots$. The $u'$-row of $\partial_k$ corresponds to paths in $\Gamma_{\!\ldots,u'}$ of length 1. Adding \smash{$\frac{-w_1'}{w_1}\!\cdot\!(u_1$}-row of $\partial_k)$ means replacing the path $u'\!\leftarrow\!v_1$ by all paths $u'\!\leftarrow\!v_1\!\to\!u_1\!\leftarrow\!v$ of length 3. Adding \smash{$\frac{-w_2'}{w_2}\!\cdot\!(u_2$}-row of $\partial_k)$ means replacing the latter by all paths $u'\!\leftarrow\!v_1\!\to\!u_1\!\leftarrow\!v_2\!\to\!u_2\!\leftarrow\!v$ of length 5. This process ends when all $v_1,v_2,\ldots$-entries of $r$ are zero (i.e. we used up all $\mathcal{M}_\leq$ to create paths). Then, taking only the $I'_k$-entries (throwing away paths that do not begin in a critical $v'$) produces $r'$.
\end{proof}
Therefore, AMT gives the same end result as Gaussian elimination with several pivots. However, with AMT we have to compute $\partial'(v)$ only for critical vertices $v$ (there may be very few of them, sometimes none at all). Also, AMT calculates the new complex recursively using (\ref{(4)}), which is much faster than row/column operations, since digraph traversal does not alter the sparse data structure.
\par As a corollary of \ref{lmm2}, we deduce that after applying AMT, the new complex $C'_\ast$ is homotopy-equivalent to $C_\ast$, an alternative proof to that of Sk\"{o}ldberg in \cite{articleSkoldbergMTFAV}.

\subsection{Examples}\label{2.3.exp} 
Several special cases of our construction \ref{2.1.formulation} have already appeared in the literature. Let us mention a few of them, along with a smaller, more illustrative instance, which also shows the importance of choosing good orders on $I_k$.
\par $\bullet$ Consider the triangulation $\Delta$ below for the real projective plane. Let $C_\ast$ be the chain complex for simplicial homology of $\Delta$. The digraph $\Gamma_{C_\ast}$ (=Hasse diagram of $\Delta$) with the lexicographic order on vertices is pictured below left, and $\Gamma_{C_\ast}$ with a different total order is shown below that. The red edges are members of the steepness matching, the circled vertices are the critical simplices, and the dotted edges / gray entries can be deleted by remark \ref{1.3.rmk}\,(c).\vspace{-2mm}
\[\includegraphics[width=0.23\textwidth]{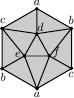}\vspace{-7mm}\]
\[\hspace{-2mm}\begin{tikzpicture}\matrix (m) [matrix of math nodes, row sep=-5pt, column sep=25pt]{
 &ab&\\
 &ac&\\
\ccircled{$a$}&ad&abd\\
 &ae&abe\\
b&a\!f&acd\\
 &\ccircled{$bc$}&ac\!f\\
c&bd&ae\!f\\
 &be&bce\\
d&\ccircled{$b_{\!}f$}&\ccircled{$bc\!f$}\\
 &cd&bdf\\
e&ce&cde\\
 &c\!f&\ccircled{$de\!f$}\\
f&de&\\
 &df&\\
 &e\!f&\\ };
\draw(m-1-2.west)--(m-3-1.east);\draw[red](m-1-2.west)--(m-5-1.east); \draw(m-2-2.west)--(m-3-1.east);\draw[red](m-2-2.west)--(m-7-1.east);
\draw(m-3-2.west)--(m-3-1.east);\draw[red](m-3-2.west)--(m-9-1.east); \draw(m-4-2.west)--(m-3-1.east);\draw[red](m-4-2.west)--(m-11-1.east);
\draw(m-5-2.west)--(m-3-1.east);\draw[red](m-5-2.west)--(m-13-1.east);\draw(m-6-2.west)--(m-5-1.east);\draw(m-6-2.west)--(m-7-1.east);
\draw[dotted](m-7-2.west)--(m-5-1.east);  \draw[dotted](m-7-2.west)--(m-9-1.east);
\draw[dotted](m-8-2.west)--(m-5-1.east);  \draw[dotted](m-8-2.west)--(m-11-1.east);
\draw(m-9-2.west)--(m-5-1.east);          \draw(m-9-2.west)--(m-13-1.east);
\draw[dotted](m-10-2.west)--(m-7-1.east); \draw[dotted](m-10-2.west)--(m-9-1.east);
\draw[dotted](m-11-2.west)--(m-7-1.east); \draw[dotted](m-11-2.west)--(m-11-1.east);
\draw[dotted](m-12-2.west)--(m-7-1.east); \draw[dotted](m-12-2.west)--(m-13-1.east);
\draw[dotted](m-13-2.west)--(m-9-1.east); \draw[dotted](m-13-2.west)--(m-11-1.east);
\draw[dotted](m-14-2.west)--(m-9-1.east); \draw[dotted](m-14-2.west)--(m-13-1.east);
\draw[dotted](m-15-2.west)--(m-11-1.east);\draw[dotted](m-15-2.west)--(m-13-1.east);
\draw[dotted](m-3-3.west)--(m-1-2.east);\draw[dotted](m-3-3.west)--(m-3-2.east);\draw[red](m-3-3.west)--(m-7-2.east);
\draw[dotted](m-4-3.west)--(m-1-2.east);\draw[dotted](m-4-3.west)--(m-4-2.east);\draw[red](m-4-3.west)--(m-8-2.east);
\draw[dotted](m-5-3.west)--(m-2-2.east);\draw[dotted](m-5-3.west)--(m-3-2.east);\draw[red](m-5-3.west)--(m-10-2.east);
\draw[dotted](m-6-3.west)--(m-2-2.east);\draw[dotted](m-6-3.west)--(m-5-2.east);\draw[red](m-6-3.west)--(m-12-2.east);
\draw[dotted](m-7-3.west)--(m-4-2.east);\draw[dotted](m-7-3.west)--(m-5-2.east);\draw[red](m-7-3.west)--(m-15-2.east);
\draw(m-8-3.west)--(m-6-2.east);\draw(m-8-3.west)--(m-8-2.east);\draw[red](m-8-3.west)--(m-11-2.east);
\draw(m-9-3.west)--(m-6-2.east);\draw(m-9-3.west)--(m-9-2.east);\draw(m-9-3.west)--(m-12-2.east);
\draw(m-10-3.west)--(m-7-2.east);\draw(m-10-3.west)--(m-9-2.east);\draw[red](m-10-3.west)--(m-14-2.east);
\draw(m-11-3.west)--(m-10-2.east);\draw(m-11-3.west)--(m-11-2.east);\draw[red](m-11-3.west)--(m-13-2.east);
\draw(m-12-3.west)--(m-13-2.east);\draw(m-12-3.west)--(m-14-2.east);\draw(m-12-3.west)--(m-15-2.east);
\end{tikzpicture}		\hspace{4mm} \raisebox{12mm}{$
R^6\!\xleftarrow{\!\!\!\left[\begin{smallmatrix}
	\mns1\!&\mns1\!&\mns1\!&\mns1\!&\mns1\!&0&\gr{0}&\gr{0}&0&\gr{0}&\gr{0}&\gr{0}&\gr{0}&\gr{0}&\gr{0}\\
	\rd{1}&0&0&0&0&\mns1\!&\gr{\mns1}\!&\gr{\mns1}\!&\mns1\!&\gr{0}&\gr{0}&\gr{0}&\gr{0}&\gr{0}&\gr{0}\\
	0&\rd{1}&0&0&0&1&\gr{0}&\gr{0}&0&\gr{\mns1}\!&\gr{\mns1}\!&\gr{\mns1}\!&\gr{0}&\gr{0}&\gr{0}\\
	0&0&\rd{1}&0&0&0&\gr{1}&\gr{0}&0&\gr{1}&\gr{0}&\gr{0}&\gr{\mns1}\!&\gr{\mns1}\!&\gr{0}\\
	0&0&0&\rd{1}&0&0&\gr{0}&\gr{1}&0&\gr{0}&\gr{1}&\gr{0}&\gr{1}&\gr{0}&\gr{\mns1}\!\\
	0&0&0&0&\rd{1}&0&\gr{0}&\gr{0}&1&\gr{0}&\gr{0}&\gr{1}&\gr{0}&\gr{1}&\gr{1}
	\end{smallmatrix}\right]\!\!\!}\! R^{15}\!\xleftarrow{\!\left[\begin{smallmatrix}
	\gr{1}&\gr{1}&\gr{0}&\gr{0}&\gr{0}&\gr{0}&\gr{0}&\gr{0}&\gr{0}&\gr{0}\\
	\gr{0}&\gr{0}&\gr{1}&\gr{1}&\gr{0}&\gr{0}&\gr{0}&\gr{0}&\gr{0}&\gr{0}\\
	\gr{\mns1}\!&\gr{0}&\gr{\mns1}\!&\gr{0}&\gr{0}&\gr{0}&\gr{0}&\gr{0}&\gr{0}&\gr{0}\\
	\gr{0}&\gr{\mns1}\!&\gr{0}&\gr{0}&\gr{1}&\gr{0}&\gr{0}&\gr{0}&\gr{0}&\gr{0}\\
	\gr{0}&\gr{0}&\gr{0}&\gr{\mns1}\!&\gr{\mns1}\!&\gr{0}&\gr{0}&\gr{0}&\gr{0}&\gr{0}\\
	0&0&0&0&0&1&1&0&0&0\\
	\rd{1}&0&0&0&0&0&0&1&0&0\\
	0&\rd{1}&0&0&0&\mns1\!&0&0&0&0\\
	0&0&0&0&0&0&\mns1\!&\mns1\!&0&0\\
	0&0&\rd{1}&0&0&0&0&0&1&0\\
	0&0&0&0&0&\rd{1}&0&0&\mns1\!&0\\
	0&0&0&\rd{1}&0&0&1&0&0&0\\
	0&0&0&0&0&0&0&0&\rd{1}&1\\
	0&0&0&0&0&0&0&\rd{1}&0&\mns1\!\\
	0&0&0&0&\rd{1}&0&0&0&0&1\\
\end{smallmatrix}\right]\!}\!R^{10}$}\vspace{-5mm}\]
\[\hspace{-2mm}\begin{tikzpicture}\matrix (m) [matrix of math nodes, row sep=-5pt, column sep=25pt]{
 &ab&\\
 &ad&\\
\ccircled{$a$}&bd&abd\\
 &ac&acd\\
b&cd&abe\\
 &ae&bce\\
d&be&cde\\
 &\ccircled{$bc$}&ac_{\!}f\\
c&ce&bc\!f\\
 &de&bd_{\!}f\\
e&a\!f&ae\!f\\
 &c\!f&\ccircled{$de\!f$}\\
f&b_{\!}f&\\
 &df&\\
 &e\!f&\\};
\draw(m-1-2.west)--(m-3-1.east);\draw[red](m-1-2.west)--(m-5-1.east);
\draw(m-2-2.west)--(m-3-1.east);\draw[red](m-2-2.west)--(m-7-1.east);
\draw[dotted](m-3-2.west)--(m-5-1.east);\draw[dotted](m-3-2.west)--(m-7-1.east);
\draw(m-4-2.west)--(m-3-1.east);\draw[red](m-4-2.west)--(m-9-1.east);
\draw[dotted](m-5-2.west)--(m-7-1.east);\draw[dotted](m-5-2.west)--(m-9-1.east);
\draw(m-6-2.west)--(m-3-1.east);\draw[red](m-6-2.west)--(m-11-1.east);
\draw[dotted](m-7-2.west)--(m-5-1.east);\draw[dotted](m-7-2.west)--(m-11-1.east);
\draw(m-8-2.west)--(m-5-1.east);\draw(m-8-2.west)--(m-9-1.east);
\draw[dotted](m-9-2.west)--(m-9-1.east);\draw[dotted](m-9-2.west)--(m-11-1.east);
\draw[dotted](m-10-2.west)--(m-7-1.east);\draw[dotted](m-10-2.west)--(m-11-1.east);
\draw(m-11-2.west)--(m-3-1.east);\draw[red](m-11-2.west)--(m-13-1.east);
\draw[dotted](m-12-2.west)--(m-9-1.east);\draw[dotted](m-12-2.west)--(m-13-1.east);
\draw[dotted](m-13-2.west)--(m-5-1.east);\draw[dotted](m-13-2.west)--(m-13-1.east);
\draw[dotted](m-14-2.west)--(m-7-1.east);\draw[dotted](m-14-2.west)--(m-13-1.east);
\draw[dotted](m-15-2.west)--(m-11-1.east);\draw[dotted](m-15-2.west)--(m-13-1.east);
\draw[dotted](m-3-3.west)--(m-1-2.east);\draw[dotted](m-3-3.west)--(m-2-2.east);\draw[red](m-3-3.west)--(m-3-2.east);
\draw[dotted](m-4-3.west)--(m-2-2.east);\draw[dotted](m-4-3.west)--(m-4-2.east);\draw[red](m-4-3.west)--(m-5-2.east);
\draw[dotted](m-5-3.west)--(m-1-2.east);\draw[dotted](m-5-3.west)--(m-6-2.east);\draw[red](m-5-3.west)--(m-7-2.east);
\draw(m-6-3.west)--(m-7-2.east);\draw(m-6-3.west)--(m-8-2.east);\draw[red](m-6-3.west)--(m-9-2.east);
\draw(m-7-3.west)--(m-5-2.east);\draw(m-7-3.west)--(m-9-2.east);\draw[red](m-7-3.west)--(m-10-2.east);
\draw[dotted](m-8-3.west)--(m-4-2.east);\draw[dotted](m-8-3.west)--(m-11-2.east);\draw[red](m-8-3.west)--(m-12-2.east);
\draw(m-9-3.west)--(m-8-2.east);\draw(m-9-3.west)--(m-12-2.east);\draw[red](m-9-3.west)--(m-13-2.east);
\draw(m-10-3.west)--(m-3-2.east);\draw(m-10-3.west)--(m-13-2.east);\draw[red](m-10-3.west)--(m-14-2.east);
\draw[dotted](m-11-3.west)--(m-6-2.east);\draw[dotted](m-11-3.west)--(m-11-2.east);\draw[red](m-11-3.west)--(m-15-2.east);
\draw(m-12-3.west)--(m-10-2.east);\draw(m-12-3.west)--(m-14-2.east);\draw(m-12-3.west)--(m-15-2.east);
\end{tikzpicture}		\hspace{4mm} \raisebox{12mm}{$
R^6\!\xleftarrow{\!\!\!\left[\begin{smallmatrix}
	\mns1\!& \mns1\!& \gr{0}& \mns1\!& \gr{0}& \mns1\!& \gr{0}& 0& \gr{0}& \gr{0}& \mns1\!& \gr{0}& \gr{0}& \gr{0}& \gr{0}\\
	\rd{1}& 0& \gr{\mns1}\!& 0& \gr{0}& 0& \gr{\mns1}\!& \mns1\!& \gr{0}& \gr{0}& 0& \gr{0}& \gr{\mns1}\!& \gr{0}& \gr{0}\\
	0& \rd{1}& \gr{1}& 0& \gr{1}& 0& \gr{0}& 0& \gr{0}& \gr{\mns1}\!& 0& \gr{0}& \gr{0}& \gr{\mns1}\!& \gr{0}\\
	0& 0& \gr{0}& \rd{1}& \gr{\mns1}\!& 0& \gr{0}& 1& \gr{\mns1}\!& \gr{0}& 0& \gr{\mns1}\!& \gr{0}& \gr{0}& \gr{0}\\
	0& 0& \gr{0}& 0& \gr{0}& \rd{1}& \gr{1}& 0& \gr{1}& \gr{1}& 0& \gr{0}& \gr{0}& \gr{0}& \gr{\mns1}\!\\
	0& 0& \gr{0}& 0& \gr{0}& 0& \gr{0}& 0& \gr{0}& \gr{0}& \rd{1}& \gr{1}& \gr{1}& \gr{1}& \gr{1}
	\end{smallmatrix}\right]\!\!\!}\! R^{15}\!\xleftarrow{\!\left[\begin{smallmatrix}
	\gr{1}& \gr{0}& \gr{1}& \gr{0}& \gr{0}& \gr{0}& \gr{0}& \gr{0}& \gr{0}& \gr{0}\\
	\gr{\mns1}\!& \gr{\mns1\!}& \gr{0}& \gr{0}& \gr{0}& \gr{0}& \gr{0}& \gr{0}& \gr{0}& \gr{0}\\
	\rd{1}& 0& 0& 0& 0& 0& 0& 1& 0& 0\\
	\gr{0}& \gr{1}& \gr{0}& \gr{0}& \gr{0}& \gr{1}& \gr{0}& \gr{0}& \gr{0}& \gr{0}\\
	0& \rd{1}& 0& 0& 1& 0& 0& 0& 0& 0\\
	\gr{0}& \gr{0}& \gr{\mns1}\!& \gr{0}& \gr{0}& \gr{0}& \gr{0}& \gr{0}& \gr{1}& \gr{0}\\
	0& 0& \rd{1}& \mns1\!& 0& 0& 0& 0& 0& 0\\
	0& 0& 0& 1& 0& 0& 1& 0& 0& 0\\
	0& 0& 0& \rd{1}& \mns1\!& 0& 0& 0& 0& 0\\
	0& 0& 0& 0& \rd{1}& 0& 0& 0& 0& 1\\
	\gr{0}& \gr{0}& \gr{0}& \gr{0}& \gr{0}& \gr{\mns1}\!& \gr{0}& \gr{0}& \gr{\mns1}\!& \gr{0}\\
	0& 0& 0& 0& 0& \rd{1}& 1& 0& 0& 0\\
	0& 0& 0& 0& 0& 0& \rd{\mns1}\!& \mns1\!& 0& 0\\
	0& 0& 0& 0& 0& 0& 0& \rd{1}& 0& \mns1\!\\
	0& 0& 0& 0& 0& 0& 0& 0& \rd{1}& 1
	\end{smallmatrix}\right]\!}\!R^{10}$}\vspace{-7mm}\]
By the first ordering, $C'_{\!\ast}\!=\! R\!\xleftarrow{\![\begin{smallmatrix}\!0&\!0\!\end{smallmatrix}]\!}\! R^2\!\xleftarrow{\!\!\left[\begin{smallmatrix}1&\mns1\\\!\mns1&\mns1\end{smallmatrix}\!\right]\!\!}\!R^2\!$, but by the second, $C'_{\!\ast}\!=\! R\!\xleftarrow{\![0]\!}\! R\!\xleftarrow{\![\mns2]\!}\!R$, with 3 critical simplices and 4 zig-zag paths, which is optimal. 
\par $\bullet$ In the example \ref{1.2.exp}, our $\mathcal{M}$ is the steepness matching.
\par $\bullet$ Sk\"{o}ldberg's matching in \cite[pp.48]{articleSkoldbergHHLAFCT} that computes the homology of Heisenberg Lie algebras is $\mathcal{M}_\leq$ w.r.t. lexicographic ordering for $x_1\!\!>\!\ldots\!>\!x_n\!>\!y_1\!\!>\!\ldots\!>\!y_n\!>\!z$.
\par $\bullet$ Forman's matching in \cite[pp.17]{articleFormanUGDMT}, which computes the homology of the simplicial complex on $E(K_n)$ of non-connected subgraphs of the full graph $K_n$, is $\mathcal{M}_\leq$ w.r.t. the lexicographic order using $12\!<\!13\!<\!23\!<\!14\!<\!24\!<\!34\!<\!15\!<\!\ldots\!<\!45\!<\!\ldots$.
\par $\bullet$ Kozlov's matching in \cite[pp.193]{bookKozlovCAT}, which computes the homology of the simplicial complex on $V(P_n)$ of all independent sets in the path graph $P_n$, is $\mathcal{M}_\leq$ w.r.t. the lexicographic order for $0\!<\!3\!<\!6\!<\!9\!<\!\ldots\!<\!2\!<\!5\!<\!8\!<\!11\!<\!\ldots\!<\!1\!<\!4\!<\!7\!<\!10\!<\!\ldots$. \hfill$\lozenge$

\subsection{Remarks}\label{2.4.rmk} Consider the following chain complexes of $\Z$-modules:\vspace{-4pt}
\[
(a)\; 0\!\leftarrow\!\Z^2\!\! \overset{\big[\!\begin{smallmatrix}3&4\\2&3\end{smallmatrix}\!\big]}{\longleftarrow} \!\!\Z^2 \!\leftarrow\!0\hspace{9mm}
(b)\; 0\!\leftarrow\!\Z^n\!\! \xleftarrow{\!\!\left[\!\!\begin{smallmatrix}
&&\hspace{-3pt}1&\hspace{-2pt}1\\
&\hspace{-6pt}\rotatebox[origin=c]{100}{$\ddots$}\hspace{-4pt} &\hspace{-12pt}\rotatebox[origin=c]{100}{$\ddots$}\hspace{-4pt}&\\[-2pt]
&\hspace{-8pt}1&\hspace{-12pt}1&\\
\,1&\hspace{-7pt}1&&\\[0pt]
\,\rd{1}&&&\\[2pt]\end{smallmatrix}\!\right]\!\!} \!\!\Z^n \!\!\leftarrow\!0\hspace{9mm}
(c)\; 0\!\leftarrow\!\Z^2 \!\!\overset{\big[\!\begin{smallmatrix}0&\!1\\\rd{1}&\!1\end{smallmatrix}\!\big]}{\longleftarrow}\!\!\Z^2 \!\leftarrow\!0\vspace{-5pt}\]
\par$\bullet$ The boundary matrices might have all entries non-invertible (so the steepness matching is empty w.r.t. any ordering), yet the chain complex can still be contractible. E.g. $(a)$ is contractible (since $\det\partial_1\!=\!1\!\in\!\Z^\times$), but  $\mathcal{M}_\leq\!=\!\emptyset$.
\par$\bullet$ The boundary matrices may have all nonzero entries invertible and be very sparse, but if $\leq$ is badly chosen, the steepness matching can still be very small. E.g. for $(b)$ our $\mathcal{M}\!=\!\{n\!\leftarrow\!1\}$ has one edge. However, if we pick the same order on rows and reverse order on columns, $\partial_1$ becomes upper-triangular and $\mathcal{M}\!=\!\{i\!\leftarrow\!i;\,i\!\in\![n]\}$.
\par$\bullet$ The inclusion in Lemma \ref{lmm1}\,(b) is not an equality. E.g. for $(c)$ the steepness pairing $\mathcal{M}\!=\!\{2\!\leftarrow\!1\}$ is not maximal. It is strictly contained in the Morse Matching $\mathcal{M}\!=\!\{2\!\leftarrow\!1,1\!\leftarrow\!2\}$, which is the steepness pairing when the rows are switched.
\par$\bullet$ The boundary matrices may have all nonzero entries invertible, but the steepness matching could be small no matter the choice of the ordering. E.g. for the complex \smash{$0\!\leftarrow\!\Z^m\!\! \overset{\partial}{\longleftarrow} \!\!\Z^n \!\leftarrow\!0$} in which every entry of $\partial$ is $\pm1$, any choice of $\leq$ gives $\mathcal{M}\!=\!\{m\!\leftarrow\!1\}$, though admittedly this matrix $\partial$ is not sparse.
\par$\bullet$ Needless to say, for most chain complexes in practice, e.g. those coming from simplicial complexes (Poincar\'{e}), (semi)groups (Eilenberg-MacLane), associative algebras (Hochschild), Lie algebras (Chevalley), knots and links (Khovanov), etc., the matrices are usually very sparse, with just entries $0$ and $\pm1$. Hence the steepness matching is very large and useful, as it kills the majority of basis elements. \hfill$\lozenge$\vspace{4mm}

\section{Choice of ordering}\label{3.ordering}
In this section, we describe how different orders on the sets of basis elements $I_k$ impact the calculation in AMT, and how to choose good orders.\\

When $\partial_k$ are given as finite matrices, the bases and orders are already determined: $I_k\!=\!\{1,\ldots,|I_k|\}$ and $\leq_k$=(usual total order on $\N)$. We can get any other order by permuting the rows and/or columns of boundary matrices.
\par Even in very sparse complexes, bad orders lead to small matchings and therefore slow computation (as seen in remark \ref{2.4.rmk}\,$(b)$), or they lead to increase in matrix densities and therefore memory overflow (as seen in the following example).
\subsection{Examples}\label{3.1.exp} Consider the 3 orderings of the same chain complex:\vspace{-3pt}
\[\hspace{-5pt}C_{\!\ast}\!:\; \; \overset{\raisebox{5mm}{$\displaystyle(a)$}}{0}\!\!\leftarrow\!\Z^m\!\! \xleftarrow{\!\!\left[\!\begin{smallmatrix} 1&\\[-0pt]\scalebox{0.8}{\rotatebox{90}{$.\hspace{0.4pt}.\hspace{0.4pt}.$}}&&\!0\!\\1\\\rd{1}&\!1\!&\!\cdots\!&\!1 \end{smallmatrix}\!\right]\!\!} \!\!\Z^n \!\!\leftarrow\!0\hspace{7mm}
\overset{\raisebox{5mm}{$\displaystyle(b)$}}{0}\!\!\leftarrow\!\Z^m\!\! \xleftarrow{\!\!\left[\!\begin{smallmatrix} \!1\!&\!1\!&\!\cdots\!&\!1\\1\\\scalebox{0.8}{\rotatebox{90}{$.\hspace{0.4pt}.\hspace{0.4pt}.$}}&&\!0\!\\\rd{1}&\\ \end{smallmatrix}\!\right]\!\!} \!\!\Z^n \!\!\leftarrow\!0\hspace{7mm}
\overset{\raisebox{5mm}{$\displaystyle(c)$}}{0}\!\!\leftarrow\!\Z^m\!\! \xleftarrow{\!\!\left[\!\begin{smallmatrix} \rd{1}\!&\!\cdots\!&\!1\!&\!1\\&&&\!1\\&\!0\!&&\!\scalebox{0.8}{\rotatebox{90}{$.\hspace{0.4pt}.\hspace{0.4pt}.$}}\\&&&\!\rd{1}\\\end{smallmatrix}\!\right]\!\!} \!\!\Z^n \!\!\leftarrow\!0\vspace{-2mm}\]
\[\hspace{-5pt}C'_{\!\ast}\!:\; 0\!\leftarrow\!\Z^{m_{\!}-_{\!}1}\!\! \xleftarrow{\!\!\left[\!\begin{smallmatrix} \mns1&\!\cdots\!&\mns1\\[-0pt]\scalebox{0.8}{\rotatebox{90}{$.\hspace{0.4pt}.\hspace{0.4pt}.$}}& &\scalebox{0.8}{\rotatebox{90}{$.\hspace{0.4pt}.\hspace{0.4pt}.$}}\\
	\mns1&\!\cdots\!&\mns1\end{smallmatrix}\!\right]\!\!} \!\!\Z^{n_{\!}-_{\!}1} \!\!\!\leftarrow\!0\hspace{5mm}
0\!\leftarrow\!\Z^{m_{\!}-_{\!}1}\!\! \xleftarrow{\!\!\left[\!\begin{smallmatrix} 1\!&\!\cdots\!&\!1\\&\phantom{0}\\&\!0\!\\\end{smallmatrix}\!\right]\!\!} \!\!\Z^{n_{\!}-_{\!}1} \!\!\!\leftarrow\!0\hspace{5mm}
0\!\leftarrow\!\Z^{m_{\!}-_{\!}2}\!\! \xleftarrow{\!\!\left[\!\begin{smallmatrix} &\!\phantom{a}\!&\\\!\phantom{a}\!&\!0\!&\!\phantom{a}\!\\&\!\phantom{a}\!& \\\end{smallmatrix}\!\right]\!\!} \!\!\Z^{n_{\!}-_{\!}2} \!\!\!\leftarrow\!0\vspace{-0mm}\]
If in $(a)$ we switch the first and last row, we get $(b)$; if we also switch the first and last column, we get $(c)$. The starting matrix density is $\frac{m+n-1}{mn}$. After computing $C'_\ast$, in $(a)$ the density increased to $1$, in $(b)$ it decreased to $\frac{1}{m\!-\!1}$, and in $(c)$ it became $0$. \hfill$\lozenge$\\

\subsection{Algorithm} To meaningfully reduce $C_\ast$, we must create as many red entries like in (\ref{(6)}), as few situations like in example \ref{3.1.exp}\,$(a)$, and as few zig-zag paths like in (\ref{(7)}) (i.e. new entries in $C'_\ast$) as possible. This is achieved in two steps:
\begin{gather}\begin{tabular}{l}For $k\!=\!1,\ldots,N$, permute the columns of $\partial_k$ (and thus rows of $\partial_{k+1}$, with \\
the same permutation), by lexicographically comparing, for every column\\
index $v$, the tuple $\big(c_0,c_1,c_2,c_3\big)(v)$. Here we used:\\
$c_0$ = 0 if the $v$-column contains a unit entry and 1 otherwise,\\
$c_1$ = position of the last nonzero entry in the $v$-column,\\
$c_2$ = density of the $v$-column,\\
$c_3$ = 0 if the last nonzero entry in the $v$-column is a unit and 1 otherwise.
\end{tabular} \tag{8} \label{(8)}\vspace{-2mm}\end{gather}
\begin{gather}\begin{tabular}{l}For $k\!=\!N\!-\!a,N\!-\!a\!-\!b,N\!-\!a\!-\!2b,N\!-\!a\!-\!3b,\ldots$ where $a\!\in\!\{0,1\}$, $b\!\in\!\{1,2\}$, \hspace{0pt}\phantom{.}\\
permute rows of $\partial_k$ (and thus columns of $\partial_{k-1}$, with the same permutation),\\
by lexicographically comparing, for every row index $u$, \\
the tuple $\big(r_0,r_1,r_2,r_3\big)(u)$. Here we used:\\
$r_0$ = 1 if the $u$-row contains a unit entry and 0 otherwise,\\
$r_1$ = position of the first nonzero entry in the $u$-row,\\
$r_2$ = -(density of the $u$-row),\\
$r_3$ = 1 if the first nonzero entry in the $u$-row is a unit and 0 otherwise.
\end{tabular} \raisetag{15mm} \tag{9} \label{(9)}\vspace{-0mm}\end{gather} 
We can use only certain $r_i$ or $c_i$ and in a different sequence. For any ordered set $I\!\subseteq\!\{0,1,2,3\}$, let (\ref{(8)})$_{I}$ and (\ref{(9)})$_{I,a,b}$ denote the procedures (\ref{(8)}) and (\ref{(9)}) in which we compare lexicographically the corresponding $c_i$ and $r_i$ for $i\!\in\!I$.
\par Performing (\ref{(9)})$_{\{1\},0,1}$, i.e. comparing only $r_1$, produces a complex, in which every matrix has the form (\ref{(7)}). Similarly, applying (\ref{(8)})$_{\{1\}}$, i.e. comparing only $c_1$, produces a complex, in which every matrix has the form (\ref{(10)}) below: white area are zeros, gray area and dashed lines are any entries, full lines are nonzeros, invertible red entries are members of $\mathcal{M}_\leq$. Further comparing the other $r_i$ (resp. $c_i$) means permuting the rows (columns) in the dark gray area of (\ref{(7)}) (resp. (\ref{(10)})).
\par As an illustration, notice that in examples \ref{2.3.exp}, the second ordering of the chain complex for the projective plane has the form (\ref{(7)}).\vspace{-1mm}
\[\includegraphics[width=0.5\textwidth]{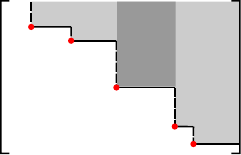} \tag{10} \label{(10)}\vspace{-2mm}\]
\par Notice that if we apply (\ref{(8)})$_{\{1\}}$ and (\ref{(9)})$_{\{1\},a,1}$, then the latter will not be reordering matrices necessarily of the form (\ref{(10)}), since it changes $\partial_k$ as well as $\partial_{k-1}$. However, if we apply (\ref{(8)})$_{\{1\}}$ and (\ref{(9)})$_{\{1\},a,2}$, then the latter always starts with the form (\ref{(10)}). In other words, we reorder the columns in all matrices, but reorder rows in only half of the matrices, to have control over the structure of the input.

\begin{Lmm}\label{lmm3} \textbf{(a)} If there exist unit entries in a matrix $\partial_k$, then after reordering it w.r.t. $(c_0,c_1,c_3)$ and $(r_0,r_1,r_3)$, the resulting matching $\mathcal{M}_\leq$ is nonempty.\\
	\textbf{(b)} Over any ring, the following process terminates eventually. While there exist unit entries in some $\partial_{k_0}$, perform (\ref{(8)})$_{\{0,1,3,2\}}$ and (\ref{(9)})$_{\{0,1,3,2\},a,2}$, $a\!=\!\big\{\begin{smallmatrix}0&\text{if }N-k_0\in2\N\\1&\text{if }N-k_0\notin2\N\end{smallmatrix}$, then reduce $C_\ast$ with AMT to a smaller $C'_\ast$ and let $C_\ast\!:=\!C'_\ast$.
\end{Lmm}\vspace{-2mm}
\begin{proof} \textbf{(a)} After reordering w.r.t. $c_0$ and $r_0$, the matrix $\partial_k$ has a $2\!\times\!2$ block form. In the $(2,1)$-block, every row and every column contains a unit entry. Further reordering w.r.t $c_1$ or $r_1$, the $(2,1)$-block obtains the form (\ref{(10)}) or (\ref{(7)}). Further reordering w.r.t. $c_3$ or $r_3$ causes the corner of the last or first step to be a unit.\\
\textbf{(b)} The choice of $a$ and $b$ assures us that we permute the columns and rows of $\partial_{k_0}$ without changing the column order again. By (a), the associated $\mathcal{M}_\leq$ is nonempty, so AMT removes at least one column and row, hence the algorithm is finite.
\end{proof}	

\subsection{Remarks} $\bullet$ If we work over a field, we use (\ref{(8)})$_{\{1,2\}}$ and (\ref{(9)})$_{\{1,2\},a,b}$. If we want to maximize the size of $\mathcal{M}_\leq$ when $\partial_\ast$ contains noninvertible entries, we use (\ref{(8)})$_{\{0,1,3,2\}}$ and (\ref{(9)})$_{\{0,1,3,2\},a,b}$. If we want to minimize the density of matrices in $C'_\ast$, we use (\ref{(8)})$_{\{0,1,2,3\}}$ and (\ref{(9)})$_{\{0,1,2,3\},a,b}$. Note however, that the choice of $a$ and $b$ can be tricky, since (\ref{(9)}) orders only half of the matrices, and may mess up the other half. Experimenting is advisable.
\par$\bullet$ In my experience, if the matrices of $C_\ast$ come from a typical homology theory (e.g. from CW complexes, groups, associative or Lie algebras, etc.), then already the default lexicographic order is often very effective. Without spending effort on reordering, $C_\ast$ is reduced a bit less, but faster. However, if $C_\ast$ consists of only one matrix (this case has important applications in \ref{5.4.SparseLinearAlgebra}), then reordering rows and columns proved vital. An experiment, in which we picked 100 random sparse matrices over $R\!=\!\Z_2$ of size $10^6\!\times\!10^6$ and density $\frac{30}{10^6}$, showed that the size of $\mathcal{M}_\leq$ after reordering was on average 6 times larger than before reordering.
\par$\bullet$ In a theoretical context, an interesting problem is to derive some characterization that determines for which orders $\leq$ is our Morse matching $\mathcal{M}_\leq$ maximal.\vspace{4mm}

\section{Noninvertible entries}\label{4.nonunits}
In this section, we describe how AMT can be used even when boundary matrices have all entries noninvertible. The obtained result is the torsion in homology.\\

Let $R$ be a PID (principal ideal domain) and $p\!\in\!R$ a prime. By the classification of finitely generated modules over PIDs, we have $H_k(C_\ast;R)\cong\bigoplus_{t\in T_k}R/(t)$ for a unique finite multiset $T_k$ of zeros and prime powers. Then over the localized ring, there holds $H_k(C_\ast;R_{(p)}\!)\cong\bigoplus_{t\in T_k\cap Rp}R_{(p)}/(t)$, which reveals $p$-torsion and free part. Computing $H_k(C_\ast;R_{(p)}\!)$ for all (relevant) primes $p$ then produces $H_k(C_\ast;R)$.
\par The only noninvertible elements in the subring $R_{(p)}\!=\!\{\frac{a}{b}\!\in\!Q(R); b\!\notin\!Rp\}$ of the field of fractions of $R$ are $R^\times_{(p)}\!=\!\{\frac{a}{b}\!\in\!Q(R); a\!\in\!Rp,b\!\notin\!Rp\}$, the multiples of $p$.\\

\subsection{Algorithm}\label{4.1.algorithm} As long as a chain complex contains a unit entry $u\!\overset{w}{\leftarrow}\!v$, there is a Morse matching $\mathcal{M}$ containing that entry (e.g. $\mathcal{M}\!=\!\{u\!\overset{w}{\leftarrow}\!v\}$ suffices).
\par Let $C_\ast^{[0]}\!:=\!C_\ast$. Computing over $R_{(p)}$, we have a sequence (\ref{(5)}), where $C_\ast^{(r)}$ contains only nonunit entries (i.e. multiples of $p$). Let $p^{a_1}$ be the largest power of $p$ that divides all nonzero entries of $C_\ast^{(r)}$ (if there are any). Let $C_\ast^{[1]}$ be the complex $C_\ast^{(r)}$ in which all entries are divided by $p^{a_1}$; then $C_\ast^{[1]}$ contains unit entries. Thus we can apply (\ref{(5)}) again, get a complex with all entries nonunits, divide it by the largest divisor $p^{a_2}$, and obtain a complex $C_\ast^{[2]}$ that contains unit entries. Doing the same again produces $C_\ast^{[3]}$ and so on. Let $m_{k,i}\!\in\!\N$ be the sum of $|\mathcal{M}_k|$, where $\mathcal{M}$ runs through all Morse matchings, used when applying (\ref{(5)}) to transform $C_\ast^{[i-1]}$ into $C_\ast^{[i]}$.\vspace{-0mm}

\begin{Lmm}\label{lmm4} The above sequence of complexes terminates at some $C_\ast^{[t]}$ in which all matrices are zero. Then the $p$-torsion of $H_k(C_\ast;R)$ is $\bigoplus_{1\leq i<t}(R/\!Rp^{a_1+\ldots+a_i}\!)^{m_{k_{\!}+_{\!}1,i}}$\!.
\end{Lmm}\vspace{-2mm}
\begin{proof} The sequence terminates, since every $e\!\in\!\mathcal{M}_k$ reduces the finite rank of $C_k$ and $C_{k_{\!}-\!1}$ by 1. Via using \ref{AMT}\,(b) repetitively, the above procedure reconstructs the SNF of $C_\ast$ over $R_{(p)}$. Indeed, by constructing $C_\ast^{[1]}\!,\ldots,C_\ast^{[t]}$, we successively change the bases of all $C_k$, such that the boundary matrices obtain the form\vspace{-1mm}
\[\begin{array}{l}
\partial_k\!=\! \mathrm{diag}(\id_{m_{k,0}},0,p^{a_1}\partial_k^{[1]}),\\
\partial_k\!=\! \mathrm{diag}(\id_{m_{k,0}},0,p^{a_1}\id_{m_{k,1}},0,p^{a_{1\!}+a_2}\partial_k^{[2]}),\\
\partial_k\!=\! \mathrm{diag}(\id_{m_{k,0}},0,p^{a_1}\id_{m_{k,1}},0,p^{a_{1\!}+a_2}\id_{m_{k,2}},0,p^{a_{1\!}+a_2+a_3}\partial_k^{[3]}),\\[-5pt]
\hspace{14pt}\vdots\\[-3pt]
\partial_k\!=\! \mathrm{diag}(\id_{m_{k,0}},0,p^{a_1}\id_{m_{k,1}},0,\ldots,p^{a_{1\!}+\ldots+a_{t\!-\!1}}\id_{m_{k,t\!-\!1}},0,p^{a_{1\!}+\ldots+a_t}0),\vspace{-1mm}
\end{array}\]
where the last zero is the trivial map $C_{k_{\!}-\!1}^{[t]}\!\leftarrow\!C_k^{[t]}$ that determines the free part.
\end{proof}

\subsection{Remark} In my experience, if the matrices of $C_\ast$ over $\Z$ have only entries $0$ and $\pm1$, and if the torsion part of $H_\ast C_\ast$ is small (relative to the ranks of $C_k$), then applying AMT produces $C'_\ast$ in which most entries are still $0$ and $\pm1$. Consequently, for many complexes  it is effective to work over $\Z$ directly, not over $\Z_{(p)}$ for all $p$.
\par For complexes $C_\ast$ with large torsion part in homology, it would help if one could efficiently compute some finite (reasonably small) set of primes $P$, such that $H_\ast C_\ast$ contains no $p$-torsion for $p\!\in\!\N\!\setminus\!P$.\vspace{2mm}

\subsection{Reconstructing SNFs}\label{4.3.algorithm} With the above algorithm, it is possible to recover the SNFs of all boundary matrices in $C_\ast$ and thus also $H_k(C_\ast;R)$. 
\par The procedure goes as follows. Over the field $Q(R)$, repeatedly apply AMT to $C_\ast$ until the resulting boundaries $\partial'_k$ are zero. If $\partial'_k$ has dimensions $|I'_{k-1}|\!\times\!|I'_k|$, denote $r_k=\sum_{0\leq i<k}(-1)^{k-1-i}(|I_i|\!-\!|I'_i|)$. Let $S_k$ be the empty multiset. Go through all primes $p$ and for each one, compute the SNF of all $\partial_k$ over $R_{(p)}$, as in the proof of \ref{lmm4}. Each time, add the obtained elementary divisors to $S_k$, then convert $S_k$ to the multiset of invariant factors $\widetilde{S_k}$. Terminate the process when $|\widetilde{S_k}|\!=\!r_k$ for all $k$. 

\begin{Lmm}\label{lmm5} In the above, the resulting $\widetilde{S_k}$ are the invariant factors of $\partial_k$ over $R$. 
\end{Lmm}\vspace{-2mm}
\begin{proof} When $C_\ast$ is transformed into $C'_\ast$, by \ref{AMT} we have $|I_k|\!=\!|I'_k|\!+\!r_k\!+\!r_{k+1}$ for all $k$, where $r_k$ is the rank of $\partial_k$ over $Q(R)$. Consequently, we can express the rank as \vspace{-1mm} \[\textstyle{r_k= \sum_{i=0}^k(-1)^{k-i}(r_i\!+\!r_{i-1})  =\sum_{i=0}^k(-1)^{k-i}(|I_{i-1}|\!-\!|I'_{i-1}|),}\vspace{-1mm}\] as in the procedure above. Since rank over $Q(R)$ equals the number of invariant factors over $R$, the algorithm terminates when all SNFs have been recovered. 
\end{proof}\vspace{4mm}

\section{Computer Implementation}
In this section, we sum up our findings from previous sections on how to compute the reduced chain complex $C'_\ast$. Then, we compare the efficiency of our basic \textsc{Mathematica} implementation of AMT to the efficiency of homology algorithms in other CASs (computer algebra systems). 

\subsection{Algorithm}\label{5.1.algorithm} We have a procedure to compute everything from \ref{AMT}.\\[-3mm]
\begin{enumerate}[label={(\alph*)},leftmargin=5mm,align=left] 
	\item[] Input: Matrices $\partial_\ast\!=\!(\partial_k)_{k=1}^N$ representing the complex $C_\ast$ over $R$.\\ 
	Let $I_{k-1}$ (resp. $I_k$) be the list of row (resp. column) indices of $\partial_k$. At start, let $f_\ast\!:=\!(f_k)_{k=0}^N$ and $g_\ast\!:=\!(g_k)_{k=0}^N$, where $f_k\!=\!g_k$ is the $|I_k|\!\times\!|I_k|$ identity matrix.\\
	\item[] While at least one of the matrices $\partial_k$ contains at least one invertible entry:
	\begin{enumerate} \setlength{\itemindent}{-0mm}
		\item[(a)] Permute the basis elements (rows of $\partial_k$ and columns of $\partial_{k-1}$ with the same permutation $\pi_k$) for all $k$, so that as many invertible entries have zeros left and below them. For instance, like in (\ref{(8)}) and (\ref{(9)}).
		\item[(b)] Compute $\mathcal{M}_\leq\!=\!(\mathcal{M}_1,\ldots,\mathcal{M}_N)$, where $\mathcal{M}_k$ is a list $\big(u_1\!\overset{w_1}{\leftarrow}\!v_1,\ldots,u_r\!\overset{w_r}{\leftarrow}\!v_r\big)$ of all invertible entries in $\partial_k$ that have zeros below in their columns and left in their rows, as defined in (\ref{(6)}). Denote $I_k^-\!=\!(u_1,\ldots,u_r)$, $I_k^+\!=\!(v_1,\ldots,v_r)$, as defined in \ref{formulation}. Compute \smash{$W_k\!=\!\big(\frac{-1}{w_1},\ldots,\frac{-1}{w_r}\big)$}, the inverses in $R$, and $I'_k\!=\!I_k\!\setminus\!(I_k^+\!\cup\!I_{k+1}^-)$. For all $k$, set all the entries of $I_k^-$-columns in $\partial_{k-\!1}$ and $I_k^+$-rows in $\partial_{k+1}$ to zero. They are not needed, by remark \ref{1.3.rmk}\,(c); this simplifies the search for paths $\gamma$ in $\Gamma_{C_\ast}^{\mathcal{M}}$. Furthermore, multiply the $I_k^-$-rows of $\partial_k$ by $W_k$. Then $\partial_\gamma$ will be computed via the simpler formula in \ref{1.3.rmk}\,(c), as weights of $e\!\in\!\mathcal{M}$ are now redundant. 
		\item[(c)] From $\partial_k,\mathcal{M}_k$ compute $\partial'_k,\, f'_k,\, g'_{k_{\!}-_{\!}1}$, for all $k$, using formulas (\ref{(4)}) and \ref{1.3.rmk}\,(b). 
		The most demanding step (searching for weighted paths $\gamma$ that give $\partial_\gamma$) can be done as follows. Set the $\mathcal{M}_k$-entries in $\partial_k$ to zero. Replace the list $\mathcal{M}_k$ (its elements are edges \smash{$u\!\overset{w}{\leftarrow}\!v$}) with a dictionary data structure (its key-value pairs are $\mathcal{M}_k(u)\!=\!v$). Replace the sparse matrix $\partial_k$ with a list-of-lists of column entries (its v-th entry is $\partial_k(v)\!=\!\big((u_1,w_1),(u_2,w_2),\ldots\big)$ of all nonzero entries in the $v$-column of $\partial_k$). The queries $\mathcal{M}_k(u)$ and $\partial_k(v)$ are completed instantaneously, since the implementation of a table and dictionary offer this. Create an empty list of entries $\partial'_k\!=\!()$. 
		\subitem To compute the matrix $\partial'_k$, go through all $v\!\in\!I'_k$. We must find all paths $\gamma$ from $v$ to all $u\!\in\!I'_{k-1}$. To do this, create a list $F\!=\!()$ of final paths and $T$ of temporary paths; initially, $T$ consists of paths with 1 edge, namely from $v$ to all its neighbors. Then recursively use $\partial_k$ and $\mathcal{M}_k$ to move between $I_k$ and $I_{k-1}$, performing a depth-first traversal. In the $i$-th step of the recursion, pick a path $\gamma\!\in\!T$ and let $u_i\!\in\!I_{k-1}$ be its endpoint: if $u_i\!\in\!I'_{k-1}$, move $\gamma$ from $T$ to $F$; otherwise, replace $\gamma$ in $T$ with all paths, obtained from $\gamma$ by adding 2 edges (one edge via $\mathcal{M}_k(u_i)\!=\!v_i$, one via $\partial_k(v_i)\!=\!\{(u_{i+1},w_{2i}),\ldots\}$). When $T$ becomes empty, return $F$. 
		It is unnecessary to store whole paths $v\!\overset{w_0}{\to}\!u_1\!\!\overset{\!1}{\to}\!\! v_1 \!\overset{w_2}{\to}\!u_2\!\!\overset{\!1}{\to}\!\! v_2\!\overset{w_4}{\to}\!\ldots\! \overset{\!1}{\to}\!\!v_i\overset{\!w_{2i}}{\to}\!\! u_{i+1}$ during the calculation, it suffices to only keep $(u_{i+1},w_0w_2\!\cdots\!w_{2i})$. When all paths starting in $v$ are found, group together results $(u,w_j)$ with the same endpoint $u$: append the entry $(u,v)\!\rightarrow\!\Sigma_jw_j\!\in\!R$ to the list $\partial'_k$. When all entries in the list $\partial'_k$ are computed, construct the sparse matrix $\partial'_k\!\in\!R^{|I_{k\!-\!1}|\times|I_k|}$ from those entries, and then set it to be the submatrix of $\partial'_k$ with $I_{k-1}'$-rows and $I'_k$-columns. 
		\subitem In a similar manner, we obtain $f'_k$ and $g'_{k_{\!}-_{\!}1}$, by finding paths (recursively or iteratively) from $v\!\in\!I'_k$ to $u\!\in\!I_k$ and from $v\!\in\!I_{k_{\!}-_{\!}1}$ to $u\!\in\!I'_{k_{\!}-_{\!}1}$, respectively. 
		\item[(d)] For all $k$, let $\partial_k\!:=\!\partial'_k$ (replace the old boundary with the new). Multiply the $I_k^-$-columns of $g'_{k-1}$ by $W_k$. Then permute the rows of $f'_k$ and columns of $g'_k$ by $\pi_k^{\!-\!1}$\!. Finally, let $f_k\!:=\!f_kf_k'$ and $g_k\!:=\!g'_kg_k$ (matrix multiplication).
	\end{enumerate}\vspace{1mm}
	\item[] Output: Matrices $\partial'_\ast\!=\!(\partial'_k)_{k=1}^N$ that represent the reduced chain complex $C'_\ast$, as well as matrices $f_\ast\!=\!(f_k)_{k=0}^N$ and $g_\ast\!=\!(g_k)_{k=0}^N$ (the h-equivalences $C_\ast\!\longleftrightarrows\!C'_\ast$, so that there holds $g_\ast\!\circ\!f_\ast\simeq\id_{C'_\ast}$ and $f_\ast\!\circ\!g_\ast\simeq\id_{C_\ast}$). 
\end{enumerate}\vspace{1mm}

\subsection{Remarks} If $R$ is a field, this process returns zero matrices $\partial'_k$. Then the width of $\partial'_k$ equals $\dim H_kC_\ast$ and columns of $f_k$ are the generators of $H_kC_\ast$, by (\ref{(5)}).
\par$\bullet$ Step (a) is optional: without it, the matching is smaller. Also, there are many ways to perform (a), so it leaves a lot of room for optimization. When matrices $\partial_k$ come from common homology theories, the default lexicographic order is often quite good for computation, and additional reordering improves the calculation a little, but spends unnecessary time. However, if we work over a PID and seek torsion, then reordering is crucial (if we perform \ref{lmm3}, then all entries in $\partial'_\ast$ are nonunits). Also, if the complex consists of only one matrix, then reordering significantly enlarges the matching and speeds up computation.
\par$\bullet$ The part of (b) about setting the columns in $\partial_{k-1}$ and rows in $\partial_{k+1}$ to zero and multiplying rows in $\partial_k,g'_{k-1}$ is optional: if we skip it, the calculation will just require more time and memory. In particular, each time a path $\gamma$ is constructed, we will additionally have to check if its endpoint lies in $I^+_{k-1}$ and compute the inverse in $R$ of a weight $w$ of an edge in $\mathcal{M}$ (possibly several times for the same edge). 
\par$\bullet$ The part of (c) and (d) about $f_\ast,g_\ast$ is optional: if we only want the isomorphism type of $H_kC_\ast$, we compute just $\partial'_\ast$; if we also want generators, we compute  $\partial'_\ast$, $f_\ast$.
\par$\bullet$ After implementing the algorithm, a good way to test its correctness is to apply it to any complex $C_\ast$ over a field, until we get $\partial'_\ast\!=\!0$. Then, check if $\partial'_{k-1}\partial'_k\!=\!0$ (i.e. $\partial'_\ast$ is a boundary), $\partial_kf_k\!=\!0$, $g_{k-1}\partial_k\!=\!0$ (i.e. $f_\ast$ and $g_\ast$ are chain maps), and $g_kf_k\!=\!\id$ (i.e. $g_\ast\!\circ\!f_\ast\simeq\id_{C'_\ast}$) for all $k$.


\subsection{Complexity}\label{5.2.complexity} As we have seen in examples \ref{2.3.exp}, without reordering the memory requirements of AMT for some complexes can be as high as if the matrices were dense. However, simply permuting the matrix rows/columns solves this problem. 
\par As seen in the motivation subsection, in general the densities converge to zero as matrices get larger. Notice that the sparser the matrices are, the larger the Morse matching is, hence the fewer critical vertices we need to deal with. Consequently, AMT becomes proportionally more effective in higher dimensions. However, it can still happen that applying AMT several times causes the matrices to fill up so much at later stages that the algorithm bogs down. In such a case, a switch to dense algorithms should be made. It also often happens that applying AMT just once reduces the complex so much that little work is needed to finish the job. It is an interesting problem to determine the expected time and memory complexity of our algorithm with reordering.

\subsection{Dynamic programming}\label{5.3.dynamicProgramming} Let us mention an alternative way to compute $C'_\ast$. 
\par By the construction of steepness matchings \ref{2.1.formulation}, the vertices in every zig-zag path $\gamma$ are decreasing w.r.t. $\leq_k$ on $I_k$. In other words, in the matrix $\partial_k$ every $\gamma$ goes from right to left and 'climbs up', like in the picture \ref{(7)}. If the paths are long, many of them might eventually join and end in the same vertex. Therefore, create an empty dictionary $D_k$. Start with the smallest $v'\!\in\!I'_k$ (left-most column index) and go onwards. Perform a depth-first search for all paths starting in $v'$. Every time a partially constructed path $\gamma$ goes through $e\!\in\!\mathcal{M}_k$, check if $D_k$ contains a key $e$. If no, compute all paths $\gamma'$ starting with $e$, add them to $D_k$ under the key $e$, then concatenate $\gamma$ with each $\gamma'$ to get all paths from $v'$. If yes, then avoid searching for paths that have already been found: concatenate $\gamma$ with all $\gamma'\!\in\!D(e)$. 
\par Whether such memoization is beneficial depends on the structure of the matrices. Sacrificing memory for $D_k$ to conserve time makes sense only when the paths are sufficiently long on average. We have not yet implemented this algorithm, so it remains to be seen whether it brings any advantages.

\subsection{Sparse linear algebra}\label{5.4.SparseLinearAlgebra} A special case of our algorithm, namely when $N\!=\!1$, offers new methods to compute the fundamental operations on sparse matrices. On sufficiently sparse matrices, these methods outperform currently built-in functions in many (perhaps all?) computer algebra systems. The results of this subsection will be treated in more detail in an upcoming article, we only give an outline here. 
\par Any matrix $\partial\!\in\!R^{m\times n}$ over a field $R$ can be regarded as a chain complex $R^m\!\overset{\partial}{\leftarrow}\!R^n$\!. From (\ref{(5)}) and the algorithm \ref{5.1.algorithm}, we obtain the sparse matrices $f_0,f_1,g_0,g_1$ below left, which induce isomorphisms on homology below right:\vspace{-2mm} \[\xymatrix{ 
R^m \ar@<2pt>[d]^{g_0}    	& R^n\ar[l]_\partial \ar@<2pt>[d]^{g_1}              & & \Coker\partial=\frac{R^m}{\Im\partial} \ar[d] & \Ker\partial=\partial^{-\!1}\!(0) \ar[d]\\
R^{m'} \ar@<2pt>[u]^{f_0}	& R^{n'}\ar[l]_{\;\;\partial'\!=0} \ar@<2pt>[u]^{f_1}& & R^{m'} \ar[u]^{f_0=g_0^{\!-\!1}}  			   & R^{n'} \ar[u]^{f_1=g_1^{\!-\!1}}}\vspace{-1mm}\]
\noindent Hence, there holds 
$\Ker g_0\!=\!\Im\partial$, $\Im f_1\!=\!\Ker\partial$, $R^m\!=\!\Im f_0\!\oplus\!\Im\partial$, $R^n\!=\!\Ker g_1\!\oplus\!\Ker\partial$. 
\par The number $m\!-\!m'\!=\!n\!-\!n'$ equals $\dim\Im\partial$, i.e. the \emph{matrix rank} of $\partial$. The columns of $f_1$ are the basis elements of the \emph{null space} $\Ker\partial$. The columns of $f_0$ are the basis elements of the \emph{quotient space} $\Coker\partial$, or equivalently, the basis elements of the \emph{complementary space} of $\Im\partial$ in $R^m$. In general, this is not an orthogonal complement, since we work over an arbitrary field $R$.
\par Let $A\!=\!\langle a_1,\ldots,a_k\rangle,B\!=\! \langle b_1,\ldots,b_l\rangle\leq R^m$ be subspaces, given by sparse generator sets (not necessarily bases), and let $\mathbf{a}$ (resp. $\mathbf{b}$) be the matrix with columns $a_1,\ldots,a_k$ (resp. $b_1,\ldots,b_l$). To compute a basis of a \emph{subspace} $A$, calculate $g_0$ for $\mathbf{a}$ (so that $A\!=\!\Ker g_0$) and $\widetilde{f_1}$ for $g_0$ (so that $\Ker g_0\!=\!\Im \widetilde{f_1}$); then the columns of $\widetilde{f_1}$ are a basis for $A$. To decide whether $A$ \emph{is contained} in $B$, calculate $g_0$ for $\mathbf{b}$ (so that $B\!=\!\Ker g_0$); then $A\!\leq\!B$ iff the product $g_0\mathbf{a}$ is the zero matrix. To compute the \emph{sum subspace} $A\!+\!B$, calculate $g_0$ for $[\mathbf{a},\mathbf{b}]\!\in\!R^{m\times(k+l)}$ (so that $A\!+\!B\!=\!\Ker g_0$) and $\widetilde{f_1}$ for $g_0$ (so that $\Ker g_0\!=\!\Im\widetilde{f_1}$); then the columns of $\widetilde{f_1}$ are a basis for $A\!+\!B$. To compute the \emph{intersection subspace} $A\!\cap\!B$, calculate $\alpha\!:=\!g_0$ for $\mathbf{a}$ and $\beta\!:=\!g_0$ for $\mathbf{b}$ (so that $A\!=\!\Ker\alpha$ and $B\!=\!\Ker\beta$), then calculate $\widetilde{f_1}$ for $\left[\begin{smallmatrix}\alpha\\\beta\end{smallmatrix}\right]$ (so that $\Ker\alpha\!\cap\!\Ker\beta\!=\!\Im\widetilde{f_1}$); the columns of $\widetilde{f_1}$ are a basis for $A\!\cap\!B$. To compute the \emph{preimage subspace} $\partial^{-\!1\!}(A)$, calculate $g_0$ for $\mathbf{a}$ (so that $A\!=\!g_0^{-1}(0)$) and $\widetilde{f_1}$ for $g_0\partial$ (so that $(g_0\partial)^{-\!1\!}(0)\!=\!\Im \widetilde{f_1}$); then the columns of $\widetilde{f_1}$ are a basis for $\partial^{-\!1\!}(A)$. When $A\!\leq\!B$, to compute the \emph{subquotient space} $B/A$, calculate $g_0$ for $\mathbf{b}$, then the columns of $f_1$ for the complex $R^{m'}\!\overset{g_0}{\longleftarrow}\! R^{m}\!\overset{\mathbf{a}}{\longleftarrow}\!R^n$ are the basis elements of $B/A$.

\subsection{Performance}\label{5.7.performance}
We test the efficiency of our \textsc{Mathematica} implementation of the algorithm \ref{5.1.algorithm} on some interesting examples. To the authors' knowledge, the homologies of some objects below (namely (a1), (a2), (b1), (b2)) have not yet been computed, or at least have not appeared in the literature.
\par All computations were carried out on a ZBook 17 G5 laptop computer with a 6-core 2.6--4.3GHz i7-8850H processor, 64GB DDR4 2666MT/s memory + 64GB swap (part of SSD disk, functioning as RAM), 512GB PCIe NVMe TLC SSD disk, and Linux Mint 19.1 operating system with \textsc{Mathematica} 11.3, \textsc{SageMath} 8.8.
\par We use the abbreviations s=second(s), m=minute(s), h=hour(s), d=day(s). For each experiment, we report the CPU and RAM performance (total time and maximum memory, needed for the calculation). \textsc{SageMath} calculations did not involve the computation of homology generators. \textsc{Mathematica} calculations did not involve reordering of rows/columns.\\[0mm]

\par\textbf{(a)} $\bullet$ Let $\Delta_n$ be the simplicial complex of all anticliques in the $n$-th hypercube graph (i.e. the independence complex of $Q_n$); it has many facets, 2 are the largest, with dimension $2^{n-1}\!-\!1$. 
Let $C_\ast$ be the Poincar\'{e} chain complex of $\Delta_n$. When $n\!=\!5$, the complex uses up 24MB of space on the disk, the largest matrix is $\partial_6\!\in\!\Z^{48\,960\times54\,304}$ with 380\,128 (nonzero) entries, its plot is below left. We computed $H_\ast(\Delta_5;\Z)$, i.e. the input are the facets of a simplicial complex, and $H_\ast C_\ast(\Delta_5;\Z)$, i.e. the input is a sparse chain complex.\vspace{-1mm}
	{\small\begin{center} \begin{tabular}{|l|l|l|l|} \hline 
		& AMT & AMT+generators & \textsc{SageMath} \\ \hline
		$H_\ast (\Delta_5;\Z)$ 		& /\!\!/ 		& /\!\!/    		 & 6h34m, 7.0GB \\
		$H_\ast C_\ast(\Delta_5;\Z)$& 3s, 30MB& 4s, 30MB & aborted after 3d, 64+10GB \\
	\hline\end{tabular}\end{center}}\vspace{-1mm}
\par $\bullet$ Let $\Delta_{m,n}$ be the simplicial complex of all non-dominating sets in the torus graph $C_m\square C_n$; it has $mn$ facets, all of dimension $mn\!-\!5$.
Let $C_\ast$ be the Poincar\'{e} chain complex of $\Delta_{m,n}$. When $m\!=\!4,n\!=\!7$, it uses up 24GB, the largest matrix is $\partial_{12}\!\in\!\Z^{22\,334\,977\times22\,118\,712}$ with 287\,543\,256 nonzero entries, its plot is below right. \vspace{-1mm}
	{\small\begin{center} \begin{tabular}{|l|l|l|l|} \hline 
		& AMT & AMT+generators & \textsc{SageMath} \\ \hline
		$H_\ast (\Delta_5;\Z)$ 		& /\!\!/ 		& /\!\!/    & crashed after 6h17m, 64+64GB \\
		$H_\ast C_\ast(\Delta_5;\Z)$& 9m, 51GB     & 14m, 64+6GB & crashed after 55m, 64+64GB  \\
	\hline\end{tabular}\end{center}}\vspace{-1mm}
\[\includegraphics[width=0.45\textwidth]{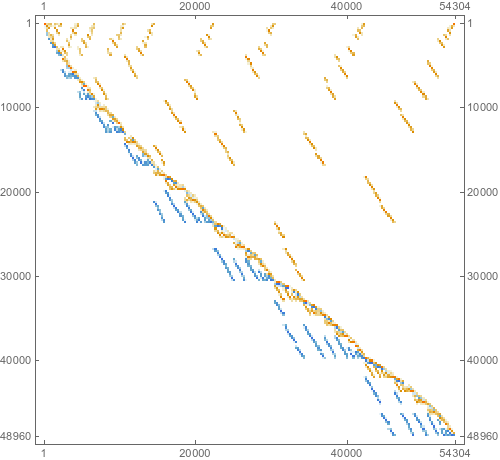}\hspace{1mm}\includegraphics[width=0.5\textwidth]{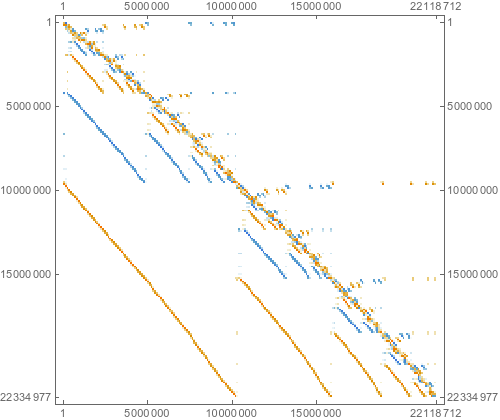}\vspace{-1mm}\]
\par $\bullet$ Let $\Delta_{m,n}$ be the $m\!\times\!n$ chessboard complex. When $m\!=\!n\!=\!8$, the dimension is $7$, the largest matrices are $\partial_5\!\in\!\Z^{376\,320\times564\,480}$ with 3\,386\,880 entries and $\partial_6\!\in\!\Z^{564\,480\times322\,560}$ with 2\,257\,920 entries; the complex uses up 122MB and contains 8\,378\,944 entries. Applying AMT 8 times spends 51m and 7GB, the new chain complex contains 2\,613\,430 entries, the largest matrices are $\partial'_5\!\in\!\Z^{42\times14\,584}$ with 36\,842 entries and $\partial'_6\!\in\!\Z^{14\,584\times8\,333}$ with 2\,576\,588 entries. The density of $\partial_6$ increased from $1\!\cdot\!10^{-6}$ to $2\!\cdot\!10^{-2}$; due to this, another application of AMT takes forever to finish, so I aborted.
A similar story happens for the matching complex (i.e. the independence complex of the line graph of the complete graph $K_n$). 
\vspace{1mm}
\par\textbf{(b)} $\bullet$ Let $A_n=\Z[S_{n,n}]$ be the category algebra over $\Z$ of the symmetric semigroup $S_{n,n}\!=\!\{\text{maps }[n]\!\to\![n]\}$ w.r.t. composition; $A_n$ has dimension $n^n$. Let $n\!=\!3$ and $C_\ast$ the first 4 matrices of the Hochschild chain complex of $A$ with coefficients in $A$. The complex uses up 882MB, its largest matrix is $\partial_4\!\in\!\Z^{474\,552\times12\,338\,352}$ with 58\,803\,294 entries, its plot is below left.\vspace{-1mm}
	{\small\begin{center} \begin{tabular}{|l|l|l|l|} \hline 
		& AMT 		& AMT+generators& \textsc{SageMath} \\ \hline
		$H_\ast C_\ast(A_3;\Z_2)$	& 34m, 15GB & 1h37m, 28GB & 4h41m, 17GB \\
		$H_\ast C_\ast(A_3;\Q)$		& 33m, 19GB & 1h33m, 32GB & crashed after 1d, 64+64GB \\
	\hline\end{tabular}\end{center}}\vspace{-1mm}
\noindent Note that with reordering, the algorithm reduced $C_\ast$ over $\Z$ to a complex with ranks 6,\,5,\,5,\,6,\,11\,881\,377 and $\partial_4$ having 100\,421 entries (all $\pm2,\pm3,\pm4,\pm6$), so an implementation of \ref{4.1.algorithm} would recover the whole $H_\ast C_\ast(A_3;\Z)$. 
\par $\bullet$ Let $A_n=\Lambda_\Z[x_1,\ldots,x_n]/(x_1\!\cdots\!x_n)$ be the exterior Stanley-Reisner algebra over $\Z$ of the simplicial $(n\!-\!2)$-sphere; $A_n$ has dimension $2^n\!-\!1$. Let $n\!=\!3$ and $C_\ast$ the first 9 matrices of the Hochschild chain complex of $A$ with coefficients in $A$. The complex uses up 2.5GB, its largest matrix is $\partial_9\!\in\!\Z^{11\,757\,312\times70\,543\,872}$ with 134\,369\,262 entries, its plot is below right. \vspace{-1mm}
	{\small\begin{center} \begin{tabular}{|l|l|l|l|} \hline 
		& AMT 	& AMT+generators& \textsc{SageMath} \\ \hline
		$H_\ast C_\ast(A_3;\Z_2)$& 2h43m, 46GB &  crashed after 3h, 64+64GB &  aborted after 4d, 22GB \\
		$H_\ast C_\ast(A_3;\Q)$	 & 2h43m, 49GB &  crashed after 3h, 64+64GB &  aborted after 4d, 39GB \\
	\hline\end{tabular}\end{center}}\vspace{-1mm}
\noindent Note that over $\Z$, the algorithm reduced $C_\ast$ to a complex with ranks 7, 20, 39, 65, 100, 151, 229, 348, 561, 60\,466\,527 and $\partial_9$ having 11\,760 entries (all between $-28$ and $25$), so an implementation of \ref{4.1.algorithm} would recover the whole $H_\ast C_\ast(A_3;\Z)$. \[\includegraphics[width=0.47\textwidth]{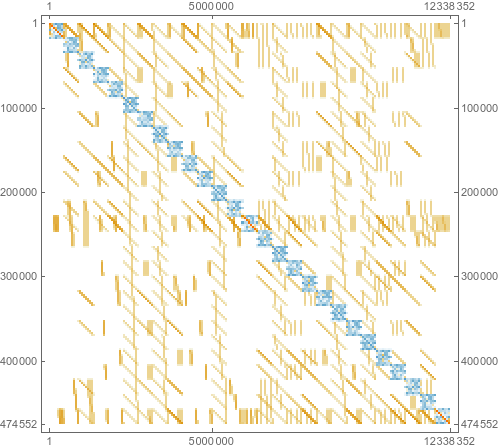}\hspace{3mm}\includegraphics[width=0.5\textwidth]{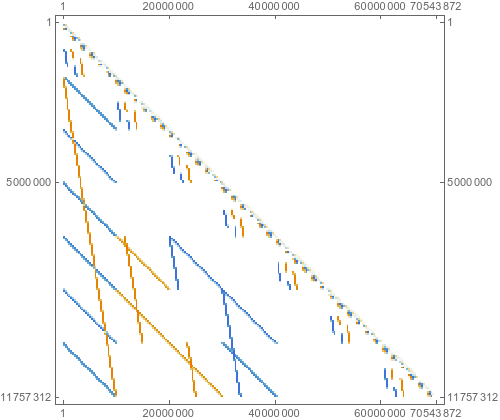}\vspace{-1mm}\]
\par\textbf{(c)} $\bullet$ Let $\frak{g}_n$ be the general linear Lie algebra over $\Z$ (it has dimension $n^2$) and $C_\ast$ its Chevalley chain complex. When $n\!=\!5$, it uses up 5.9GB, the largest matrix is $\partial_{13}\!\in\!\Z^{5\,200\,300\times5\,200\,300}$ with 66\,339\,260 entries, its plot is below left.\vspace{-0mm}
	{\small\begin{center} \begin{tabular}{|l|l|l|l|} \hline 
		& AMT & AMT+generators & \textsc{SageMath} \\ \hline
		$H_\ast C_\ast(\frak{g}_5;\Z_2)$& 35m, 14GB & 55m, 14GB & aborted after 30d, 26GB \\ 
		$H_\ast C_\ast(\frak{g}_5;\Q)$	& 15m, 14GB & 24m, 14GB & crashed after 8d, 64+18GB \\
	\hline\end{tabular}\end{center}}\vspace{-0mm}
\par $\bullet$ Let $\frak{g}_n$ be the Heisenberg Lie algebra over $\Z$ (it has dimension $2n\!+\!1$) and $C_\ast$ its Chevalley chain complex. 
When $n\!=\!13$, the complex uses up 4.5GB, the largest matrix is $\partial_{14}\!\in\! \Z^{20\,058\,300\times20\,058\,300}$ with 35\,154\,028 entries, its plot is below right.\vspace{-0mm}
	{\small\begin{center} \begin{tabular}{|l|l|l|l|l|} \hline 
		& AMT 		  & AMT+generators & \textsc{SageMath} \\ \hline
		$H_\ast C_\ast(\frak{g}_{13};\Z_2)$ & 34m,   18GB & 1h9m,  21GB & aborted after 10d, 24GB \\
		$H_\ast C_\ast(\frak{g}_{13};\Q)$   & 1h47m, 18GB & 4h19m, 21GB & aborted after 10d, 56GB \\
	\hline\end{tabular}\end{center}}\vspace{-0mm}
\[\includegraphics[width=0.47\textwidth]{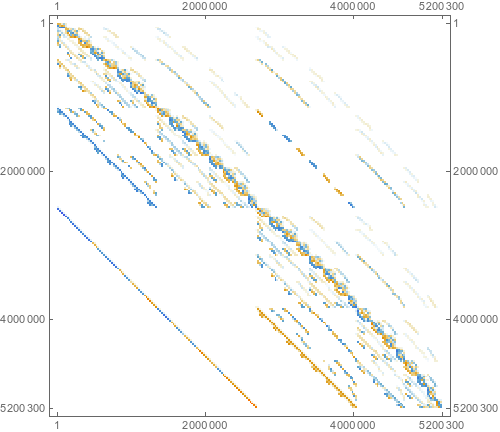}~\hspace{3mm}\includegraphics[width=0.49\textwidth]{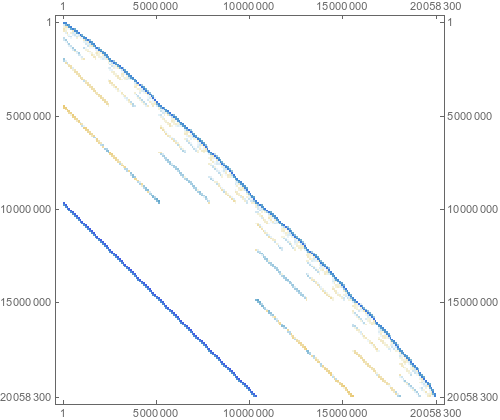}\]
\vspace{4mm}

\section{Conclusion} As we have seen in \ref{5.7.performance}, a rudimentary implementation of our new algorithm \ref{5.1.algorithm} succeeds in computing the homology of very large complexes on just a laptop and with no reordering of the columns/rows of matrices. In most cases, it worked very fast and with small memory requirements (for instance, it computed the homology over $\Z$ of the independence complex of the $8,3$-Kneser graph, with dimension $20$ and largest matrix of size 5 million $\times$ 5 million). However, in some instances (usually when the ranks of $C_k$ are large, $N$ is small, and $H_\ast C_\ast$ contains a lot of torsion), in the new complex the density of matrices increases to the amount that AMT slows down too much to be efficient. I suspect that for such cases, a very specific reordering would remedy the problem, though how to choose such orders remains an open problem.\vspace{4mm}

\section{Acknowledgment}
I wish to thank Ale\v{s} Vavpeti\v{c} for the idea of considering the set of edges with the most upward tilt. Its generality really sparked my interest in this topic.
\par This research was supported by the Slovenian Research Agency program P1-0292 and grants J1-7025 and J1-8131.

\end{document}